 \documentclass[3p,times,preprint]{elsarticle} 

\usepackage{amssymb}
\usepackage{amsmath}
\usepackage{empheq}
\usepackage{orcidlink}

\usepackage{endnotes}
\let\footnote=\endnote 

\journal{Finance Research Open}

\begin{document}

\begin{frontmatter}  



\title{From Ponzi Schemes to Benign Investment Dynamics:\\
modelling Collapse, Stability, and a Path to Sustainability}

\author{Bernhard R. Parodi\, \,\orcidlink{0009-0000-9094-0407}} 
\affiliation{
	organization={GIBZ Gewerblich-industrielles Bildungszentrum~Zug},
          addressline={Baarerstrasse~100}, 
          postcode={CH-6301},
          city={Zug},
          country={Switzerland\\ 20 February 2026}
}
\ead{bernhard.parodi@gibz.ch}

\begin{abstract}
The population and capital dynamics of three stylized investment systems are
mathematically described using discrete-time difference equations with closed-form
solutions. The models share a common capital budget equation but differ in their
demographic laws, which are geometric, quasi-logistic, or epidemiologic
(SIR-based). The quasi-logistic model is designed as an analytically tractable
non-Ponzi investment system: it generalizes the geometric model (and, in the
limit of a constant growth rate, reproduces classical Ponzi dynamics) while
closely mirroring the behaviour of an SIR-based model with decreasing effective
growth. In all cases, promised returns are modeled as fixed per-period payouts
on initial investment with principal repaid upon exit, so that aggregate
liabilities depend only on the current number of active investors. Within this unified
framework, classical Ponzi schemes arise as special cases that inevitably
collapse, while suitable parameter choices in the quasi-logistic and SIR-based
versions generate finite-horizon, legally benign “no-Ponzi game” investment
schemes with analytically transparent conditions for collapse, stability, and
sustained operation. 
\end{abstract}

\begin{keyword}
Dynamic systems \sep System dynamics \sep Financial fraud \sep Scheme sustainability \sep Investment strategy \sep Risk management \sep Economic education
\JEL C6 \sep C73 \sep G1 \sep G3
\end{keyword}

\end{frontmatter} 

\twocolumn  


\section{Introduction}
\label{Sect1}

\paragraph{Summary}
Three analytic models of investment schemes embedding modified Ponzi-like payout structures are developed and contrasted. These systems feature competing cash inflows and outflows (independent of interest accrual), which---absent sufficient initial capital, market returns, demographic scale, or optimal lock-up periods---lead to bankruptcy; otherwise, they may sustain long-term activity. A couple of models admit resilient equilibria with positive outcomes, illuminating parameter thresholds that distinguish collapse from viability and offering insights into recurrent financial instruments akin to bonds or banking reserves. Highlights:
\begin{itemize}
\item Three discrete-time models (and one related conti\-nuous-time model) with closed-form solutions for capital dynamics.
\item Investor growth follows geometric, quasi-logistic, or SIR-based demography.
\item System profitability or Ponzi-like collapse based on parameters (initial capital, market exposure, lock-up times, demographic rate).
\item Identifies critical transitions between unsustainable and viable regimes.
\item Proposes a novel temporary instrument viable under controlled conditions.
\end{itemize}

\paragraph{Background}
Classical Ponzi schemes, lacking genuine market exposure or external revenue streams, offer unrealistically high returns financed exclusively through new investors’ contributions. Such mechanisms are inherently unsustainable and collapse once inflows weaken. Two principal payout configurations are commonly observed: (i) fixed periodic coupons, in which investors receive constant payments per period with principal remaining constant until redemption — the framework adop\-ted in this study; and (ii) compounding account structures, which display smoothly increasing balances through systematic reinvestment, typically accompanied by implausibly low volatility. Despite structural differences, both rely on persistent excess inflows over outflows in the absence of productive investment and often promote reinvestment to defer liquidity pressure.

Historical examples include Adele Spitzeder, Sarah Howe, and Baldomera Larra in 1870s Europe, Carlo Ponzi in the 1920 United States, and Bernard Madoff in 2008. These cases are referenced in literary works such as Dickens’s Little Dorrit and Galdós’s Cánovas \cite{Rid21} and analyzed by Sprengers \cite{Spr20} and Frankel \cite{Fra12}. Frankel \cite{Fra12} provides a psychological and historical interpretation of Ponzi schemes as recurrent socio‑econo\-mic phenomena, while Sprengers \cite{Spr20} conducts an extensive empirical investigation covering more than one thousand U.S. cases, including cryptocurrency‑related variants. Further quantitative perspectives are discussed in \emph{Introduction to the Theories and Varieties of Modern Crime in Financial Markets} \cite{Fru16}, which examines statistical characteristics of Ponzi and related fraudulent investment systems using time‑series resampling methods. Complementary detection tests based on the Sharpe ratio and Benford’s law are proposed in \cite{Zhu11} and \cite{Win21}.

The present models extend and generalize the classical coupon‑based framework into a unified dynamic capital system. Unlike all related papers cited below which focus mostly on collapse dynamics and on alleged sustainability, we provide a unified framework that reveals parameter-dependent outcomes across lawful and Ponzi-like regimes. Our formulation integrates initial capital endowments, modest market returns, finite demographic expansion, and lock‑up constraints within a single mathematical structure. All models share a common capital‑budget identity, combining compounding assets at rate $i$ with non‑compounding payout obligations linked to initial deposits. Ponzi‑type collapse emer\-ges endogenously as a limiting regime under adverse parameter configurations, while favorable regimes generate stable, bond‑like investment processes. This generalized dynamical representation formalizes the continuum between unsustainable Ponzi mechanisms and viable pooled‑income systems, thereby providing a theoretical foundation for characterizing investment dynamics as continuous outcomes of a single parametric system rather than distinct categorical phenomena. Within this broader theoretical setting, it is useful to briefly review prior research on Ponzi mechanisms and related investment systems to position the present model in the existing literature.

\paragraph{Literature review} 
This paper focuses exclusively on investment systems superordinate to classic Ponzi schemes, in particular on their underlying mathematical mechanisms and on \emph{processes by which unscrupulous schemes may be reformed to achieve benign outcomes}. In parti\-cular, a mathematical framework is provided that captures both the typological distinction made in \cite{Spr20} between intentional and unintentional Ponzi schemes\footnote{Intentional Ponzi schemes originate as deliberate fabrications designed to attract victims, whereas unintentional Ponzi schemes result from business failures. The research by \cite{Spr20}, which utilized a dataset based on U.S. federal agency documents spanning nearly six decades, found that approximately 90\% of Ponzi schemes were intentional.}---inevitably leading to disastrous outcomes--- and modified Ponzi schemes or Ponzi scheme-like systems with the possibility of benign results. Before proceeding, we provide a brief overview of related technical studies on classical Ponzi schemes, encompassing both coupon‑\linebreak based and compounding payout structures. In a recent review of unmodified classical Ponzi schemes \cite{BP25}, the theo\-retical section selectively summarizes several mathe\-matical formulations for modelling fatal capital evolution in both continuous‑ and discrete‑time frameworks. The present analysis extends and generalizes these approaches. 

For \emph{continuous-time} Ponzi models (e.g., \cite{Cla09}, \cite{Peng21}, \cite{Per21},\ \cite{BP25}) that are, moreover, analytically tractable in closed-form  ---excluding those with complicated integral equations that ultimately need to be tackled numeri\-cally---  see for example \cite{Artz09}, \cite{Kos17}, \cite{Parlar25}, and \cite{Par14} (where the latter incorporates  a lock-up period and a stochastic parameter to be set $\sigma=0$).   In \cite{Fru16} it is stated that a "Ponzi scheme is destined to fail from the outset but the timing and the amplitude of the collapse depends on many factors, like the rate of funds inflow, the withdrawal rate, and the volatility of the investment, if any." In particular, lock-up or holding periods regulate withdrawals and are a significant and often underestimated element in Ponzi schemes. In \cite{Peng21} it is observed that "similar to hedge funds, mandating such a period helps the Ponzi scheme to survive its fledgling stage and escape an early collapse". As will be highlighted in this paper, in a (non-classical or modified) system with revenues derived from the capital market and with a declining population growth rate, the proper selection of a lock-up period — among other refinements — may prevent a collapse. A novel analyti\-cal continuous-time model that recasts these features and straightforwardly follows from one of our discrete-time models is presented as supplementary material in  \ref{AppB}. 

In contrast to continuous-time approaches, this paper mainly analyzes models for Ponzi schemes in \emph{discrete time} (e.g., \cite{Bha03}, \cite{Abd09}, \cite{Car11}, \cite{Zhu11}, \cite{Cun14}, \cite{Peng21}, \cite{Win21}, \cite{Bar25}, \cite{BP25}), focusing moreover on deterministic frameworks that yield analytically tractable, closed-form solutions. These solutions employ either elementary functions (\cite{Par13}, \cite{Par17}, \cite{BP25}, \cite{Parlar25} Eq. 7) or summation formulas (\cite{Ati21} Eq. 3.4, \cite{Parlar25} Eq. 3). Some models introduce extra parameters for realism or refinement, such as the 'stylistic' model in \cite{BP25}, which involves periodic payouts to the promoter, or the model in \cite{Par13}, which incorporates multiple investor payout methods and a lock-up period. The models presented here combine both elementary and summing approaches and include several of these refinements. As supplementary material and relying on our formalism, explicit solutions of the purely recursive models presented in \cite{BP25} and \cite{Win21} — the latter implementing the auxiliary fictitious-fund construction introduced in \cite{Zhu11}  —  as well as of models with sum-based solutions as in \cite{Ati21} and \cite{Parlar25} are provided in the Notes.

To address demographic dynamics, three approaches are compared: geometric growth, quasi-logistic growth (introduced here as a novel element), and epidemiological growth represented by a variant of the SIR model with an exact solution \cite{Lem25}, along with an insightful variant. This comparative procedure reveals some similarities and differences among these related models.

\paragraph{Scope}
Classical Ponzi schemes, in their unmodified\linebreak form, are unsustainable in the long run and lead to one of two well‑known outcomes: either the assets diminish to zero, or the capital stock grows until the investor pool is exhausted, rendering obligations unfulfillable. A third scenario involves the premature appropriation of capital by the promoter, which is not examined further here. In this paper, the baseline framework generalizes these classical configurations and reinterprets their behavior within a broader class of investment models that share a common capital‑budget identity but differ in their demographic dynamics. Classical schemes then arise as a special parameter regime with a specific payout structure, characterized here by non‑compounding liabilities rather than by compounding on accumulated balances. The analysis focuses on modified schemes that yield a wider economic interpretation — a family of transitory investment games, some of which are legally benign.

Three types of participation-growth models are considered: (i) the geometric model, which reproduces the familiar exponential growth setting and includes the classical Ponzi case when the growth rate $n_t=n$ is constant; (ii) the quasi‑logistic model, which extends this by introducing a finite pool of potential investors and a time‑dependent, sigmoidally decreasing growth rate $n_t$, generating empirically realistic humped participation patterns and declining effective growth similar to SIR‑type contagion models; and (iii) the non‑standard SIR model, which offers a third, epidemiologically motivated demographic law that parallels the quasi‑logistic dynamics.

The paper’s main focus is the quasi‑logistic specification, interpreted as a non‑Ponzi investment system analytically situated between the geometric and SIR‑based cases. It preserves the tractability of the geometric model while incorporating a diminishing growth rate and finite investor pool that capture SIR‑like behavior. Within this unified framework, classical Ponzi schemes, fragile systems, and benign finite‑horizon pooled‑income products all emerge as parameter‑dependent outcomes of a single underlying budget equation and payout structure.

\paragraph{Structure}
The paper is organized as follows. Section \ref{Sect2} presents the population dynamics of three growth models in detail: geometric growth (Section \ref{Sect2p1}), quasi‑\-logistic growth (Section \ref{Sect2p2}), and SIR‑type growth (Section \ref{Sect2p3}). Section \ref{Sect3} then applies these demographic models to capital dynamics in investment systems superordinate to Ponzi schemes, beginning with an intuitive overview framed as “traffic‑light” scenarios (Section \ref{Sect3p2}). Capital evolution is analyzed separately for geometric demography (Section \ref{Sect3p3}), quasi‑\-logistic demography (Section \ref{Sect3p4p1}), and epidemiological demography (Sections \ref{Sect3p4p2} – \ref{Sect3p4p3}). For the latter two models, the critical selection of parameters determines whether the outcome represents a collapse (a Ponzi‑like failure) or a benign, sustainable investment process. By constructing a coherent sequence of such benign systems in Section \ref{Sect3p4p4}, the paper demonstrates the feasibility of a novel investment framework. Section \ref{Sect4} concludes with a summary and key insights. (An explicit table of contents is provided in \ref{AppC}.)


\section{Population dynamics}
\label{Sect2}

In order to mathematically describe the time-depen\-dent sizes of groups within a total population of investors it's enough to rely on the following basic quantities (independent or dependent variables as well as parameters):
\begin{itemize}
\item[$t$:] Discrete time, measured in periods (each of length $\Delta t=1$, e.g., months or years), starting with $t = 0$; the end of period $t$ is the beginning of period $t+1$. Hence $t$ acts both as a lower-case integer index counting the periods and as an ordinary time variable. 
\item[$T$:] Duration of an investment or  lock-up time or maturity; put differently, time-span of membership of an investor in the scheme: after $T$ periods an investor quits the system and withdraws the original investment $I_0$. The very first time this happens within a system is at time $t = T$.
\item[$N_{0}$:] The number of very first investors at time $t = 0$.
\item[$N_{t}$:] The number of all current investors at time $t\ge 0$. 
\item[$\Delta N_{t}^{in}$:] the number of new investors joining the system at time $t$. 
\item[$\Delta N_{t}^{out}$:] the number of investors quitting the system at time $t$. 
\end{itemize}
The iterative calculations for the different group sizes are as follows: The number of all investors at time $t$ is the sum of the active investors $N_{t-1}$ and of the new investors $\Delta N_t^{in}$ entering at time $t$, minus the exiting investors $\Delta N_{t}^{out}$:
\begin{equation}
\boxed{
	N_{t} = N_{t - 1}+\Delta N_{t}^{in}-\Delta N_{t}^{out}. \hspace{1cm} (t \geq 1)
} 
\label{recrelNt}
\end{equation}
In this paper, we juxtapose three dynamical situations: one with geometric growth, one with quasi-logistic growth, and one with growth patterns governed by three variants of the SIR model. Accordingly, $\Delta N_{t}^{\text{in}}$ and $\Delta N_{t}^{\text{out}}$ are specified for each case. In all but one case (namely, the classical or standard SIR model), the recursive representation admits an explicit analytical solution.


\subsection{Geometric growth}
\label{Sect2p1}

In the geometric growth model, calculation of the number of newly entering investors is based, firstly, on the following  plausibility assumption: the growing influx of investors is caused by at least some of the enrolled investors who motivate and acquire new investors to participate in the system. Therefore, starting with initial value $\Delta N_{0}^{in} = N_{0}$, the number of new investors is proportional to the number of all active investors (for times $t < T$) and slightly different later on:
\begin{eqnarray}
	\Delta N_{t}^{in} &=& n \, N_{t -1} \hspace{2cm} (1 \le t < T)\label{DNin0}\\
                                     &=& n \, N_{0}{(1 + n)}^{t - 1}\hspace{1cm}(t\ge 1)\label{DNin}
\end{eqnarray}
By design, for times $t \ge T$ the number of new investors is assumed to continue to exponentially grow according to rule \eqref{DNin} (instead of rule  \eqref{DNin0}), despite of the fact that some investors already have left the system. Some comments are appropriate:

\setlength{\leftmargini}{0.4cm}
\begin{itemize}
\item
The constant value \(n\) is called the geometric growth rate.
\item
Rule \eqref{DNin0} implies that the number of \emph{all} investors grows exponentially with growth factor \((1 + n)\) during the first interval of periods (with \(\Delta N_{t}^{\text{out}} = 0\)). Using equation \eqref{recrelNt}, we obtain the recursion
\begin{equation}\label{geomNtrec}
N_{t} = (1 + n) N_{t - 1}
\end{equation}
which yields \(N_{t} = N_{0} (1 + n)^{t}\).
\item
The same exponential growth for the number of \emph{new} investors is imposed by design to hold for all times; hence, rule \eqref{DNin} is characteristic of the model presented here. This switch of rules at time \(T\) reflects the unbroken popularity of the scheme and its continued attraction to new investors, and it considerably simplifies the analytical treatment. (This remains true in the next section on quasi-logistic growth.)
\end{itemize}
\noindent Finally, the number of exiting investors $\Delta N_{t}^{out}$ is equal to the number of investors who entered the system $T$ periods earlier and are now leaving it again, i.e.,
\begin{equation}\label{GDNout}
\Delta N_{t}^{out} = \Delta N_{t - T}^{in}. \hspace{1cm} (t \geq T)
\end{equation}

In \ref{AppendixTable1} the iterative change in the number of participants  is formally retraced in table form, based on equations \eqref{recrelNt} to \eqref{GDNout}. The following simple formulae emerge:
\begin{equation}\label{GNt}
\boxed{  N_t = N_0(1+n)^t
  \begin{cases}
    1 & 0 \le t < T, \\
    \left(1 - (1+n)^{-T}\right) & t \ge T,
  \end{cases} }
\end{equation}

\begin{align}
\Delta N_t^{\text{in}} &=
  \begin{cases}
    N_0 & t = 0, \\
    N_0(1+n)^{t-1} n & t \ge 1,\label{GNtin}
  \end{cases} \\
\Delta N_t^{\text{out}} &=
  \begin{cases}
    0 & 0 \le t < T, \\
    N_0 & t = T, \\
    N_0(1+n)^{t-T-1} n & t > T.\label{GNtout}
  \end{cases}
\end{align}

It follows that the sum of all investors who have left so far is equal to the difference between the current number of participants and the one that existed $T$ periods earlier, i.e., $\sum_{k=0}^t N_k=N_t-N_{t-T}$.


\subsection{Quasi-logistic growth}
\label{Sect2p2}

\begin{figure}[!t]
\centering
\includegraphics[width=.9\linewidth]{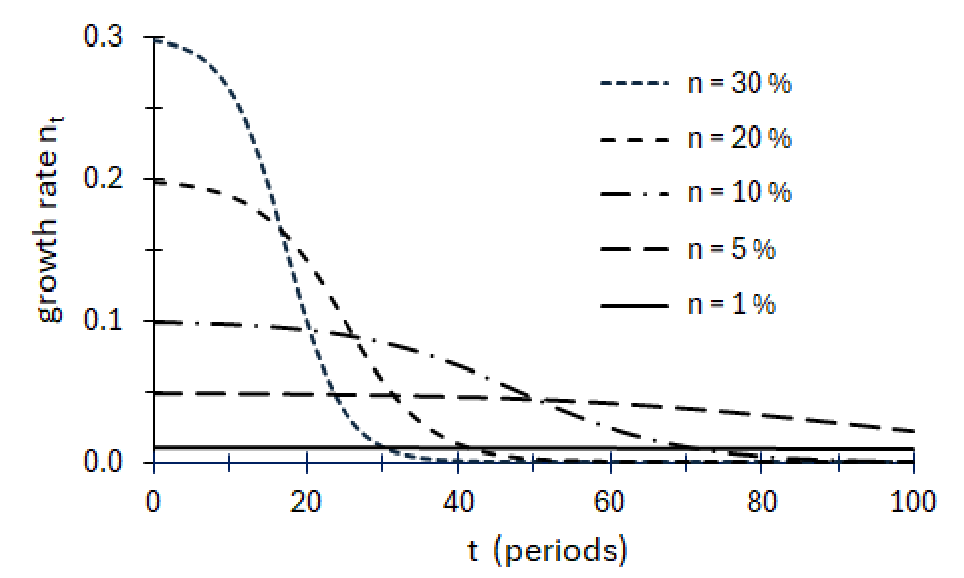}
\includegraphics[width=.9\linewidth]{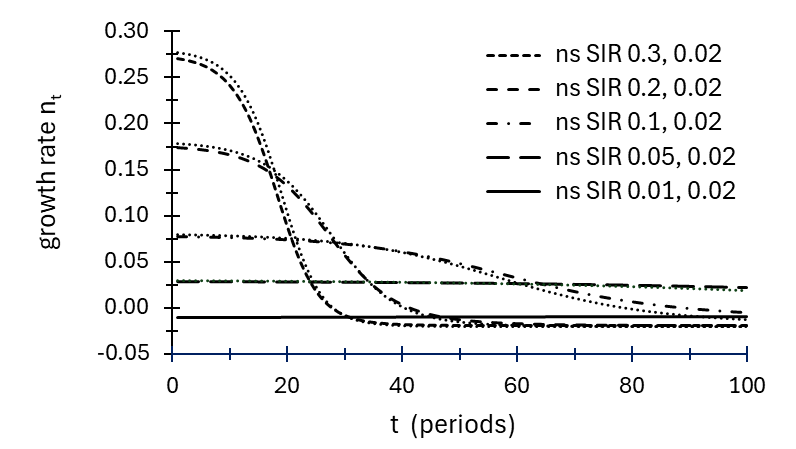}
\caption{Growth rates $n_t$ for the quasi-logistic growth model and for a couple of SIR-models. (For better visualization, only the interpolation lines connecting the discrete data points will be shown in this figure and in all other figures.) \emph{Top panel:} Sigmoidally decreasing growth rates $n_t$ for the quasi-logistic model (according to equation \ref{ntQL}). Adopting a fraction $N_0/N=0.01$, the examples shown start with initial values $n$ as indicated in the legend, pass turning points at times $t_{TP}$ (equation \ref{tTP}) and at heights $\frac{1}{2}n$, and finally approach zero. The lower the value of $n$, the larger $t_{TP}$, and hence the slower the dropping of the growth rate. \emph{Bottom panel:} Decreasing effective growth rates for two types of SIR-models, adopting the same parameter values $\beta=0.3, ..., 0.01$ and $\gamma=0.02$: the non-standard SIR-model (equation \ref{ntnsSIR}, dashed to solid lines) and, for comparison, the standard SIR-model (equation \ref{ntSIR}, dotted lines). For the purpose of illustration, equal values $\beta=n$ are chosen, hence in the case of additionally setting $\gamma=0$ the curves would coincide with those in the top panel.}
\label{Fig1ab}
\end{figure}

Realistically, the population size is bounded upwards and has a supremum $N<\infty$. This feature is exhibited by the logistic function that will therefore become a part of the design of our model. In addition, the plausibility assumption regarding new entrants (formally expressed with equation \eqref{DNin0}) is modified here in such a way that the number of new investors at time $t$ grows both proportionally to the number $N_{t-1}$ of investors still present and at a \emph{decreasing} growth rate $n_t$:
\begin{eqnarray}
\Delta N_t^{in}&=&n_t N_{t-1}   \hspace{1cm} (t\ge1),\label{DNinQL}\\
n_t&=&n\,\frac{1}{1+\frac{N_0}{N-N_0} (1+n)^t }\nonumber \\
&=&n\frac{1-\frac{N_0}{N}}{1+\frac{N_0}{N} [(1+n)^t-1] }. 
\label{ntQL}
\end{eqnarray} 

The growth rate $n_t$ decreases sigmoidally from an initial value $n_0=n (1-N_0/N)$ and asymptotically approaches zero ($\lim_{t\rightarrow \infty}n_t=0$). The form of the time-dependent (or, equivalently, period-dependent) factor is motivated such that the number of investors, $N_t$, will be given by an expression analogous to the logistic function (see below). This generates what here is referred to as quasi-logistic growth\footnote{\label{fn2}The logistic function for population growth in continuous time $t$ obeys the differential equation $dN(t)/dt=m N(t)(1-N(t)/N)$  and admits the closed-form solution $N(t)=N/[1+(N-N_0)/N\cdot e^{-mt}]$, with supremum $N$. As for the notation, the power of Euler's number $e^{m}$ is replaced here by $(1+n)$. However, the logistic function is \emph{not} a solution to the differential equation when transformed to discrete times according to replacing $N(t)/dt$ by $(N_{t+1}-N_t)/\Delta t$ (with $\Delta t=1$). Also, at later stages, our extended model predicts that the initially sigmoidal trajectory will transition into a humped profile. Therefore the designation \emph{quasi-logistic} is more accurate here.}. Figure \ref{Fig1ab} provides an illustration\footnote{\label{fn3}Introducing a trigger variable $\nu$ such that  $n_t=n\,[1+\nu\frac{N_0}{N-N_0} (1+n)^t ]^{-1}$ would allow for further controlling the steepness of the curve and hence for horizonally shifting the turning point; in the present paper a value of $\nu=1$ is set for good.}.  The intermittent turning point has a value $n_{t_{TP}}=\frac{1}{2}n$ and occurs after $t_{TP}$ periods\footnote{\label{fn4}The turning point is calculated by means of setting the second derivative to zero. Dealing with discrete time steps, we adopt the central finite difference scheme up to second order as an approximation, i.e., we demand $(n_{t+\Delta t}-2n_t+n_{t-\Delta t})/(\Delta t)^2=0$ in the limit $\Delta t\rightarrow 1$. This yields the relation $(N_0/N)(1+n)^{t_{TP}}=1-N_0/N$ which leads to both $t_{TP}$ and  $n_{t_{TP}}$.}, given by
\begin{equation}
t_{TP}=\frac{\log (\frac{N}{N_0}-1)}{ \log (1+n)}.\label{tTP}
\end{equation}
Entrant rule \eqref{DNinQL} and growth rate \eqref{ntQL} cause the number of all investors to grow sigmoidally with growth factor $(1+n_t)$ per period. During periods $0\le t<T$, equation \eqref{recrelNt} becomes 
\begin{equation}
N_t=(1+n_t)N_{t-1},
\end{equation}
which leads to $N_t\propto n_t N_0 (1+n)^t$; see Table \ref{AppendixTable2} for the formal iteration and Figure \ref{Fig1ab} for an illustration. For an infinitely large pool of members ($N\rightarrow \infty$), the formalism reduces to the case of geometric growth with constant growth rate $n_t=n$ that was discussed previously in Section \ref{Sect2p1}.

\begin{figure}[h]
\centering
\includegraphics[width=.95\linewidth]{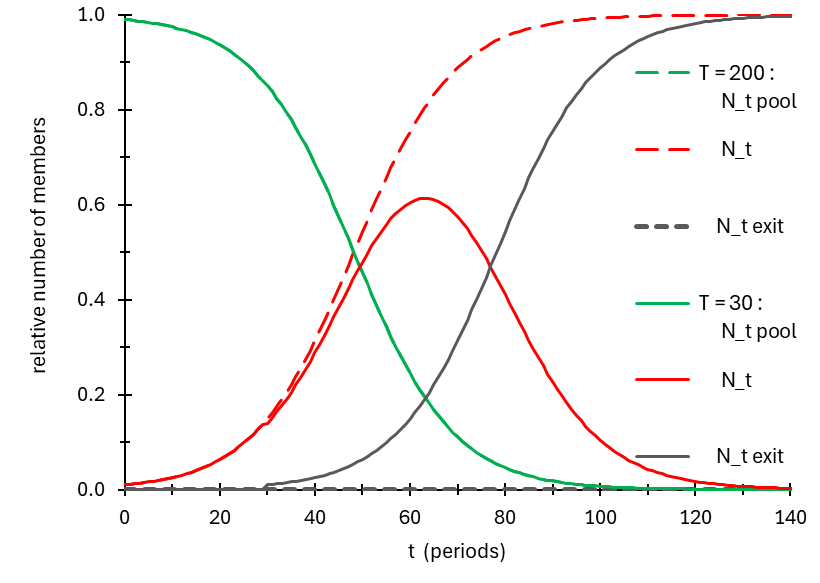}
\caption{Demographic development for the quasi-logistic model. Developments are shown under different lock-up period conditions: without exits from the system (lock-up period $T = 200$, representing $T\rightarrow \infty$, dashed lines) and with exits after $T$ = 30 periods of participation (solid lines). Sizes are normalized by $N=N_{t\,pool}+N_{t}+N_{t\,exit}$, where $N_{t\,pool}$, $N_{t}$, and $N_{t\,exit}=\sum_{k=0}^t \Delta N_k^{out}$ are the respective numbers of potential, current, and former investors in the system at time $t$.
}
\label{Fig2}
\end{figure}

The recursive equations for the number of all investors (equation \ref{recrelNt}) and for the number of exiting investors (equation \ref{GDNout}) remain formally unchanged, but the number of new investors is modified according to Equations \eqref{DNinQL} and \eqref{ntQL}. Table \ref{AppendixTable2} thus yields the following explicit calculation formulas: 
\begin{equation}\label{QLNt}
\boxed{
N_t =
	\begin{cases} 
	N_t^{t<T} 				& \hspace{1cm}  0\le t<T\\ 
	N_{t}^{t<T} - N_{t-T}^{t<T} 	& \hspace{1cm}  t \ge T
	\end{cases}
}
\end{equation}
\begin{align}
\Delta N_t^{in}&=
	\begin{cases} N_0 		& \hspace{1.8cm}  t=0 \\
	N_{t-1}^{t<T}n_t^{}		& \hspace{1.8cm}  t\ge 1\label{QLDNtin2}
	\end{cases}	\\
\Delta N_t^{out}&=
	\begin{cases} 
	0 								&    \hspace{0cm}     0\le t<T  \\
	\Delta N_{t-T}^{in} =N_{t-T-1}^{t<T}n_{t-T}^{}	&   \hspace{0cm}      t\ge T \label{QLDNtout}
	\end{cases}
\end{align}
Herein, the explicit solution for the number of currently involved investors at times $t < T$, $N_t^{t<T}$, comes into play at all times and is given by\footnote{\emph{Proof:} As exhibited in Table \ref{AppendixTable2}, for $t < T$ the number $N_t$ comprises a product over growth factors that are worked out as follows:
\small
\begin{eqnarray*}
N_t^{t < T} &=& N_0\prod_{k = 1}^{t}\left( 1 + n_k \right) \\
&=& N_0\prod_{k = 1}^{t}
\left( 1 + n \, \frac{1}{1 + \frac{N_0}{N - N_0}(1 + n)^k} \right) \\
&=& N_0\prod_{k = 1}^{t}
\left( \frac{1 + n + \frac{N_0}{N - N_0}(1 + n)^k}{1 + \frac{N_0}{N - N_0}(1 + n)^k} \right)\\ 
&=& N_0\prod_{k = 1}^{t}{(1 + n) 
\left( \frac{1 + \frac{N_0}{N - N_0}(1 + n)^{k - 1}}{1 + \frac{N_0}{N - N_0}(1 + n)^k} \right)}\\
&=& N_0(1 + n)^t \frac{1 + \frac{N_0}{N - N_0}}{1 + \frac{N_0}{N - N_0}(1 + n)^1} \cdot
\frac{1 + \frac{N_0}{N - N_0}(1 + n)^1}{1 + \frac{N_0}{N - N_0}(1 + n)^2} \cdot \\
& &\ldots \cdot \frac{1 + \frac{N_0}{N - N_0}(1 + n)^{t - 2}}{1 + \frac{N_0}{N - N_0}(1 + n)^{t-1}}
\cdot \frac{1 + \frac{N_0}{N - N_0}(1 + n)^{t - 1}}{1 + \frac{N_0}{N - N_0}(1 + n)^{t}} \\
&=& N_0(1 + n)^t \cdot \frac{1 + \frac{N_0}{N - N_0}}{1 + \frac{N_0}{N - N_0}(1 + n)^t} \\
&=& N_0(1 + n)^t \cdot \frac{1}{1 + \frac{N_0}{N}\left\lbrack (1 + n)^t - 1 \right\rbrack} \\
&=& \frac{N}{1+\frac{N-N_0}{N_0}(1+n)^{-t}}.
\hspace{0.7cm}\square
\end{eqnarray*}
\label{fn5}
Thus, the quasi-logistic sequence in the last line is merely a compact rewriting of its product representation in the first line.
\normalsize
}

\begin{equation}\label{QLNttsT}
\boxed{
N_t^{t<T} = N_0(1+n)^t\frac{1}{1-\frac{N_0}{N} +\frac{N_0}{N} (1+n)^{t}}. 
}
\end{equation}
\begin{equation}\label{QLNtlogfct}
\hspace{0cm}
                 = \frac{N}{1+\frac{N-N_0}{N_0}(1+n)^{-t}}.
\end{equation}
Expression \eqref{QLNtlogfct} is formally equivalent to the logistic function in continuous time, but here is set in action at discrete time points indexed by integers $t$, where each index denotes $t$ periods, each of duration $\Delta t=1$ ($t\equiv t\Delta t$, see Note \#\ref{fn3}). Accordingly, its values reside along a sigmoidal curve, where $\lim_{t\rightarrow \infty}N_t^{t<T} =N$ provides a supremum (red dashed line in Figure \ref{Fig2}). One may address $N_t^{t<T}$ and $N_t $ as \emph{quasi-logistic sequence} and as \emph{T-modified quasi-logistic sequence}, respectively.

In the limit of an infinitely large pool of potential investors
$N \to \infty$, the time-varying growth rate $n_t$ converges to the constant
initial value $n$, and the quasi-logistic sequence reduces to the geometric
growth law with $N_t = N_0 (1+n)^t$ discussed in Section~\ref{Sect2p1}. Formally, equation \eqref{QLNt} reduces to equation \eqref{GNt}. In this sense, the quasi-logistic model strictly embeds the geometric model as a
special case, while for finite $N$ its decreasing $n_t$ and humped $N_t$ closely parallel the effective growth rates and infection paths generated by the SIR-based models in Section~\ref{Sect2p3}.

In Figure \ref{Fig2}, the graph of the normalized number of participants $N_t/N$ shows a sigmoidal increase until the total population is reached (dashed line) in case of an unlimited investment time-span ($T \rightarrow \infty$). For a limited duration of an investment (with maturity $T < \infty$) the overall shape of the graph is dominated by a prominent hump or hill (solid line). Its maximum height ($N_{t_{Npeak}}$) or its position ($t_{Npeak}$) are determined from the condition $\Delta N_{t} = N_{t} - N_{t - 1} = 0$, with the formal result
\begin{equation}
t_{Npeak} = \frac{T+1}{2}+\frac{\ln \left( \frac{N}{N_{0}} - 1 \right)}{\ln{(1 + n)}}.
\label{QLtNpeak}
\end{equation}
For example, in Figure \ref{Fig2} $\left( t_{Npeak};N_{t_{Npeak}}/N_{total} \right) = (63.7; $ $0.614)$ (with $N/N_0=100$ and $n=10\%$). The demographic dynamics given by equations \eqref{QLNt} to \eqref{QLNttsT} serves as the basis for calculating the capital dynamics of a corresponding modified Ponzi scheme in Section \ref{Sect3p4p1}.


\subsection{Epidemiological growth (SIR model)}\label{Sect2p3}

The SIR model operates with an effectively decreasing growth rate, too. It is a mathe\-matical model in epidemiology that describes the spread of an infectious disease in a finite population. In this article, the discretized-time model serves as the basis. In the simplest version discussed here, the population is divided into the following three groups (or compartments):
\begin{itemize}
\item[$S_t$] (Susceptible): Number of people susceptible to infection; in a Ponzi scheme, these correspond to the potential investors.
\item[$I_t$] (Infected): Number of people infected and contagious; these represent the current investors in the Ponzi scheme (corresponding to the previous size $N_t$).
\item[$R_t$] (Regenerated): Number of individuals who have either recovered and are immune or have died and are redeemed; in a Ponzi scheme, these are all former members who have left the system in the meantime.
\end{itemize}
The size $N$ of the whole population remains constant in time $t \in \mathbb{N}_0$, 
\begin{equation}
N = S_0 + I_0+ R_0 = S_t + I_t + R_t.\label{Nconst}
\end{equation}

\subsubsection{Standard variant with an iterative solution}\label{Sect2p3p1}

In the standard approach without vital dynamics (i.e., birth and death rates are taken to be zero), the temporal change in group sizes is modeled by the following three coupled difference equations:
\begin{eqnarray}
S_t &=& S_{t -1} - \beta\frac{S_{t -1}}{N}I_{t -1}\label{SIRa}\\
I_t &=& I_{t -1}\  + \beta\frac{S_{t -1}}{N}I_{t -1} - \gamma I_{t -1}\label{SIRb1}\\
    &=& \left( 1 + n^{SIR}_t\right)I_{t -1}\label{SIRb2}\\
R_t &=& R_{t -1} + \gamma I_{t -1},\label{SIRc}
\end{eqnarray}
with initial values $S_0$, $I_0$, and $R_0$, and maintaining condition \eqref{Nconst} (by which $R_t$ actually bcomes redundant). Here, $\beta$ is the infection or transmission rate, $S_{t -1}I_{t -1}/N$ represents the probable fraction of newly infected persons (estimated by means of the contact probability between a susceptible and an infected person, $(S_{t -1}/N) \cdot (I_{t -1}/N)$, multiplied by the whole population $N$), and $\gamma$ is the recovery rate. The time step has a length of 1 period. In a notation analogous to the previous models, the two contributions to the change in the number of infected people are
\begin{equation}
\Delta I_t^{in}=\beta\frac{S_{t -1}}{N}\,I_{t -1}, \hspace{1cm} \Delta I_t^{out}=\gamma \, I_{t -1}.
\end{equation}
Accordingly, equation \eqref{SIRb1} can be written in the form of equation \eqref{SIRb2} by means of an effective growth rate 
\begin{equation}
n^{SIR}_t=\beta\frac{S_{t -1}}{N} - \gamma. \label{ntSIR}
\end{equation}
Because typically $S_{t -1}$ is steadily \emph{decreasing with time}, so does $n^{SIR}_t$, and as soon as $S_{t-1}=(\gamma/\beta) N$ this rate even becomes negative; thereafter, the number of infected people decreases again. In the long run, $n^{SIR}_t$ appraoches a lower limit  value $-\gamma<0$. The above standard SIR system of equations for discrete times $t$ must be solved iteratively. Figure \ref{Fig1ab} (bottom panel) illustrates the temporal evolution of the effective growth rate of infected individuals. The graphs resemble those obtained for the quasi‑logistic model (upper panel). As is well known (and not shown here), the number of infected individuals $I_t$ exhibits the characteristic hump‑shaped trajectory: the epidemic expands until it reaches its maximum spread, after which the number of infected individuals declines.

A \emph{refined} standard variant of the SIR model without vital dynamics accounts for the fact that the individuals who become recovered at time \(t\) are precisely those who were infected \(T_0\) periods earlier, whose number is given by \(I_{t-T_0}\). (In later sections --- in a capitalistic setting --- \(T_0\) is interpreted as a lock-up period for an investment.) In this refined framework, the above system of difference equations is modified by replacing \(\gamma I_{t-1}\) with
\[
\gamma\, I_{t-1-T_0}\,\mathbb{U}_{t-1-T_0},
\]
where the unit step function is defined as
\[
\mathbb{U}_x =
\begin{cases}
0, & x \le 0,\\
1, & x > 0.
\end{cases}
\]
This auxiliary piecewise function ensures the delayed onset of recoveries at the time \(t = T_0\).

In \cite{Peng21} and \cite{BP25}, where a similar model extension is proposed, the recovery rate \(\gamma\) is additionally made effectively time-dependent by replacing it with \(\gamma S_{t-T_0}/N\).\footnote{The authors formulate their model in continuous time, with the delay equal to the lock-up period (i.e., \(T_0 = \gamma^{-1}\)), and solve the resulting system of delayed differential equations numerically.} Solving such a system proceeds recursively; this direction will not be pursued further here.

\subsubsection{Non-standard variant with an exact solution}\label{Sect2p3p2}

Closely related to the standard SIR model without a lock-up time, an exact, non-negative solution is known for a variant, specifically when the difference equations \eqref{SIRa} to \eqref{SIRc} for the temporal change of group sizes are modified as follows \cite{Lem25}:
\begin{eqnarray}
S_{t} &=& S_{t - 1} - \beta\frac{S_{t}}{N - R_{t - 1}}I_{t - 1}\label{nsSIRa}\\
I_{t} &=& I_{t - 1}\  + \beta\frac{S_{t}}{N - R_{t - 1}}I_{t - 1} - \gamma I_{t} \label{nsSIRb1}\\
       &=& \left( 1 + n^{nsSIR}_t\right)I_{t -1}\label{nsSIRb2}\\
R_{t} &=& R_{t - 1} + \gamma I_{t},\label{nsSIRc}
\end{eqnarray}
for $t\ge 1$, with inital values $S_0$, $I_0$, $R_0$, and again with a time step of length 1. Compared to the standard SIR model, the factor $S_{t -1}/N$ is replaced by $S_{t}/(N - R_{t -1})$ (i.e., the former pool size from one period earlier, $S_{t-1}$, is replaced by the present size $S_t$, and it is normalized by means of only the number of remaining affected people instead of the whole population size) and accordingly the term $\gamma I_{t -1}$ was replaced by $\gamma I_t$. There obviously remains a close kinship with the original SIR-model. The temporal changes in the number of infected people are now
\begin{equation}\label{DItnsSIR1}
\Delta I_t^{in}=\beta\frac{S_{t}}{N-R_{t-1}}\,I_{t-1}, \hspace{1cm} \Delta I_t^{out}=\gamma \, I_{t}.
\end{equation}
Analogous to the quasi-logistic model with $N_{t} = (1 + n_t)N_{t - 1}$ and to equation \eqref{SIRb2} for the standard SIR model, recursion relation \eqref{nsSIRb1} can similarly be written in the form of equation \eqref{nsSIRb2}, with a decreasing \emph{effective growth rate} 
\begin{equation}
n^{nsSIR}_{t} = \frac{1}{1 + \gamma}\left(  \beta\frac{S_{t}}{N - R_{t - 1}} -\gamma \right). \label{ntnsSIR}
\end{equation}
This rate turns negative when $S_{t}=(\gamma/\beta)(N-R_{t-1})$ (with $\gamma>0$) and asymptotically approaches $-\gamma/(1+\gamma)$ for large values of $t$; a graphical illustration is given in Figure \ref{Fig1ab} (bottom panel). Given initial values $S_0$, $I_0$, and $R_0$, the \emph{explicit solution} of system \eqref{nsSIRa}$-$\eqref{nsSIRc} for $t \geq 1$ is
\begin{equation}\label{nsSIR1}
S_t = S_{0}\,p_{t},\hspace{0.5cm}
I_t = I_{0}\left( \frac{1 + \beta}{1 + \gamma} \right)^{t}p_{t},\hspace{0.5cm}
R_t = N - S_t - I_t,
\end{equation}
wherein the decisive product function
\begin{equation}
p_{t} = \prod_{k = 1}^{t}\frac{1 + \frac{S_{0\ }}{I_{0\ }}\left( \frac{1 + \gamma}{1 + \beta} \right)^{k - 1}}{1 + \beta + \frac{S_{0\ }}{I_{0\ }}\left( \frac{1 + \gamma}{1 + \beta} \right)^{k - 1}}
\label{pt}
\end{equation}
is the major prolific ingredient. In \cite{Lem25}, the proof is performed by means of complete induction. Typically, $\beta>\gamma$, hence the product function provides decreasing values, typically starting with a value very close to 1 and approximating zero for large $t$. The special case with $\gamma=0$ will become relevant in the next section. As illustrated in Figure \ref{Fig2}, the graphs of the two SIR models (the original one and its nonstandard variant) show a rather similar course. The maximum number $I_{t_{Ipeak}}$ of infected people occurs at time $t_{Ipeak}$ when condition $\Delta I_{t} = I_{t} - I_{t - 1} = 0$ holds; this delivers
\begin{equation}
t_{Ipeak}^{nsSIR} = \log\left( \frac{\gamma(1 + \gamma)}{\beta-\gamma}\,\frac{I_{0}}{S_{0}} \right) \,\,/ \,\,
\log\left( \frac{1 + \gamma}{1 + \beta} \right),
\end{equation}
requiring $\beta$ $>$ $\gamma$ $>$ 0 and $0<\frac{\gamma(1+\gamma)}{\beta - \gamma}\frac{I_0}{S_0} <1$ for a peak to occur.

The exact solution of the nonstandard SIR model, using optimized parameter values, provides a reasonable approximation to the iterative solution of the standard SIR model.

The nonstandard SIR model with a population dynamics as described by equations \eqref{nsSIR1}f. serves as the basis for calculating the capital dynamics of a corresponding modified Ponzi scheme in Section \ref{Sect3p4p2}. 

\subsubsection{Non-standard variant with delayed initial recovery}\label{Sect2p3p3}

The difference equations for the non-standard SIR model presented in the previous subsection are now modified according to subdividing the discrete timeline into two consecutive intervals of periods. In the first interval with periods $1 < t\le T_0$, no regenerations are happening that would reduce the number of currently contagious people, i.e. the recovery rate is set $\gamma=0$ in equations \eqref{nsSIRa}$-$\eqref{nsSIRc}. Regenerations according to setting $\gamma >0$ only start with period $T_0+1$, where $T_0$ defines the delay time for some initial recovery. (This parameter was already introduced in Section \ref{Sect2p3p1}). For this second, open-ended interval, the initial values are $S_{T_0}$, $I_{T_0}$, and $R_{T_0}$. Relying on the results for the unmodified nonstandard SIR model in the previous subsection, the demography here develops analogously according to the composite solution
\begin{equation}\label{nsSIR2}
S_t =
	\begin{cases} 
	S_{0}\,p_{t}^{<}\\
	S_{T_0}\,p_{t}^{>},
	\end{cases}
\boxed{
I_t =
	\begin{cases} 
	I_{0}\left( 1 + \beta \right)^{t}p_{t}^{<}  & \hspace{0.2cm}  1\le t\le T_0  \\
	I_{T_0}\left( \frac{1 + \beta}{1 + \gamma} \right)^{t-T_0}p_{t}^{>} & \hspace{0.2cm}   t \ge T_0+1,
\end{cases}
}
\end{equation}
and $R_t=N - S_t - I_t$, and wherein the corresponding decisive product functions now read\footnote{For $t\le T_0$ when $ \gamma=0$, writing out the factors for $p_t^{<}$ provides a simplification very similar to the one provided in the proof for the compact notation of $N_t^{t<T}$ in section 2.2 (cf. Note \# 5).}
\begin{eqnarray}
p_{t}^{<} &=& \prod_{k = 1}^{t}\frac{   1 + \frac{S_{0\ }}{I_{0\ }}\left( \frac{1                }{1 + \beta} \right)^{k - 1}   }
{   1 + \beta + \frac{S_{0\ }}{I_{0\ }}\left( \frac{1                 }{1 + \beta} \right)^{k - 1}   } \\
&=&\frac{1+\frac{I_0}{S_0}}{1+\frac{I_0}{S_0}(1+\beta)^t}
\label{pt1}\\
p_{t}^{>} &=& \prod_{k = T_0+1}^{t}\frac{   1 + \frac{S_{T_0}}{I_{T_0}}\left(  \frac{1 + \gamma}{1 + \beta} \right)^{k -T_0-1}   }
{   1 + \beta + \frac{S_{T_0\ }}{I_{T_0\ }}\left( \frac{1 + \gamma}{1 + \beta} \right)^{k - T_0-1}   }.
\label{pt2}
\end{eqnarray}
The less-than ($<$) and more-than symbols ($>$) are to emphasize the respective belonging to the couple of intervals. For $T_0=0$ regeneration is allowed for from the very beginning, thus by defining $p_0^{<}=1$ the first-interval solution provides the initial values and the second-inter\-val solution recovers the solution for the unmodified nonstandard SIR model (given by equations \eqref{nsSIR1}f., hence with $p_t^{>}=p_t$). The temporal changes can be calculated according to equation \eqref{DItnsSIR1}, with the solution above properly inserted.

During the first phase when $\gamma=0$, the effective growth rate in $I_t=(1+n_t^{nsSIR,<})I_{t-1}$ becomes\footnote{$n_t^{nsSIR,<}=I_t/I_{t-1}-1=(1+\beta)p_t^{<}/p_{t-1}^{<}-1=[1+\beta + I_0/S_0(1+\beta)^t]/[1 + I_0/S_0(1+\beta)^t]-1=\beta/[1 + I_0/S_0(1+\beta)^t]$.}
\begin{equation}
n_t^{nsSIR,<}=\beta\,\frac{1}{1+\frac{I_0}{S_0}(1+\beta)^t},\label{ntnsSIR1}
\end{equation}
and the number of infected people can be written as
\begin{eqnarray}
I_t^{<}&=&I_0(1+\beta)^t\,\frac{1+\frac{I_0}{S_0}}{1+\frac{I_0}{S_0}(1+\beta)^t}\\
&\approx& I_0(1+\beta)^t\,\frac{1}{1-\frac{I_0}{S_0}+\frac{I_0}{S_0}(1+\beta)^t}. \label{ItnsSIR}
\end{eqnarray}
Herein the approximation up to linear order in $I_0/S_0$ ($\approx N_0/N \ll 1$) yields expressions that are formally identical to the effective growth rate and to the number of participants in the quasi-logistic model before withdrawals (i.e., comparing equations \ref{ntnsSIR} and \ref{ntQL} and as well as \ref{ItnsSIR} and \ref{QLNttsT}). With respect to Ponzi systems, this will produce a very similar outcome  for $t<T$ and particularly for the case of perpetual participation ($T\rightarrow \infty$). More precisely, the composite solution described here (equations \ref{nsSIR2} ff.) will be used in section \ref{Sect3p4p3} as the demographic prerequisite for setting up yet another alleged or modified Ponzi system. With respect to an investment scheme, \(I_{t}\) describes the number of current investors, i.e., $N_{t} \equiv I_{t}$, and the contagion rate serves as the initial growth rate, i.e., $n=\beta$. 

We note that, unlike the model with quasi-logistic growth — where each investor exits the system after a fixed investment time span of $T$ periods (i.e., $\Delta N_t^{out}\propto N_{t-1-T}$) —  the present nonstandard SIR model assumes a different departure mechanism. Here, the group of recovered individuals leaving the "contagious" compartment (or, by analogy, investors exiting the scheme) after time 
$t=T_0$ is assumed to be proportional to the current number of infected (or actively invested) individuals. That is, $\Delta I_t^{out}\propto I_{t}$, or equivalently, $\Delta N_t^{out}\propto N_{t}$.\footnote{For the inquiry of a nonstandard model that implements not only an initial delay period but also a general lock-up time ---i.e., one in which $\Delta I_t^{out}\propto I_{t-T_0}$ for $t>T_0$--- an explicit solution is no longer available. In this case, a recursive approach would be required, which lies beyond the analytical scope of the present paper.}


\section{Capital dynamics}
\label{Sect3}

Having introduced some selected population dynamics in Section \ref{Sect2}, it now becomes embedded within a capitalistic setting. 

\subsection{Common budget equation and liability structure}\label{Sect3p1}

Setting the stage involves introducing the following new model parameters and variables related to capital evolution in general:
\begin{itemize}
\item[$I_{0}$:] one-time investment of a (new) investor when entering the scheme; this investment is profitable and acts as an interest-bearing deposit until its withdra\-wal after an ubiquitous investment time-span of $T$ periods when the investor quits the system;
\item[$r$:] rate of return on investment per period.
\item[$P$:] Periodic profit of an investor, i.e., the fixed coupon paid to every investor at the end of each period:
\begin{equation}
P = r I_{0}.
\end{equation}
The promised payout stream to an individual investor who enters at time \(t_{\mathrm{in}}\) and exits at \(t_{\mathrm{out}} = t_{\mathrm{in}} + T\) consists of
\[
\underbrace{r I_0, \; r I_0, \; \dots, \; r I_0}_{T \text{ periods}},
\quad \text{followed by} \quad I_0,
\]
with total profit \(r I_0 T\) at maturity. Alternative payout modalities that rely, for example, on contractual compounding of accumulated balances are not investigated here.\footnote{Another common payout design promises investors a large future payout based on compounded ``returns'' credited to the account: an investor contributing \(I_0\) at time \(t_{\mathrm{in}}\) is promised a lump sum at withdrawal time \(t_{\mathrm{out}} = t_{\mathrm{in}} + T\) equal to $L = I_0 (1+r)^T$, where \(r\) is the fixed per-period rate of return and \(T\) is the number of compounding periods chosen at the investor's discretion. An investor leaving the system after $T$ periods realizes a total profit of $P_T = L - I_0 = I_0 \bigl((1+r)^T - 1\bigr)$. For $r \ll 1$, a first-order approximation yields $P_T \approx r I_0 T$, which is similar in form to the total profit discussed in the main text for the coupon-based model. However, the corresponding budget equations would differ substantially---affecting the scheme's liquidity and stability---and the evolution of outstanding balances would become a different critical issue. A full modelling of such a compounding scheme lies beyond the scope of the present investigation.}
\item[$K_0^{pro}$:] Initial promoter (operator) capital contribution.
\item[$K_0$:] Total initial capital at $t=0$: 
  \begin{equation}
  K_0 = K_0^{pro} + I_0 N_0 \label{eq:K0}.
  \end{equation}
\item[$i$:] Effective mean capital market return on system assets (exogenous, post-inflation, pre-taxes, net of periodic promoter payouts). Alternatively, some value $i<0$ may directly be interpreted as a periodic payout rate to promoters.
\item[$\Delta K_t$:] Period-$t$ or net capital growth ($t\ge 1$): market returns plus new deposits minus promised payouts to prior investors minus redemptions:
  \begin{equation}
  \Delta K_t = i K_{t-1} + \Delta N_t^{in} I_0 - N_{t-1} P - \Delta N_t^{out} I_0
  \label{eq:DeltaK}
  \end{equation}
  where $\Delta N_t^{in} I_0$ and $\Delta N_t^{out} I_0$ denote aggregate subscriptions and redemptions.
Aggregate payout obligations at time $t$ therefore depend only on the current number of active
investors, $N_{t-1}$ and $N_t^{\mathrm{out}}$, and are thus given by $N_{t-1} r I_0$ (coupon payouts) and $N_t^{\mathrm{out}} I_0$ (principal repayments), respectively. Promoter payouts would introduce an additional outflow term but are absorbed here into the effective rate $i$; see \cite{Par13} for models that explicitely incorporate constant or variable promoter distributions. 
\item[$K_t$:] Accumulated capital at $t\ge 1$: $K_t = K_{0} + \sum_{s=1}^t \Delta K_s$; this is the target asset.
\end{itemize}

While obeying different demographic evolution, all three models considered in this paper share the same basic capital budget equation $K_t=K_{t-1}+\Delta K_t$ or 
\begin{equation}\label{Ktrec}
\boxed{K_t = (1+i) K_{t-1}+ \left( N_{t}^{\mathrm{in}} - N_{t}^{\mathrm{out}}- r N_{t-1}  \right) I_0},
\end{equation}
for $t \ge 1$, with initial value $K_0$. The factor $(1+i)$ implies that the 
pooled capital $K_t$ compounds at the effective 
market rate $i$ on the asset side. By contrast, promised payouts on the
liability side are non-compounding: each investor receives a fixed coupon
$r I_0$ per period on the initial deposit $I_0$, and the principal $I_0$ is
repaid once upon exit after $T$ periods. Aggregate obligations at time $t$ 
therefore depend on the current number of active investors, $N_{t-1}$ and
$N_t^{\mathrm{out}}$, but not on accumulated balances (cf. Note \# 10). 
In this paper, ``compounding'' always refers to the asset side via $i$, while the
liabilities are deliberately modeled via $r$ as non-compounding payout obligations.

This structure is identical in the geometric, quasi-logistic, and SIR-based
settings; only the demographic laws governing \(N_t\), \(N_t^{\mathrm{in}}\),
and \(N_t^{\mathrm{out}}\) differ. The triple of models considered in
Sections~\ref{Sect3p3}--\ref{Sect3p4} is therefore best viewed as a unified
family of pooled-income investment games built on a common budget equation,
with different demographic dynamics. The quasi-logistic model in particular is
designed as an analytically tractable non-Ponzi investment system intermediate
between the geometric and SIR-based specifications. 

This modelling choice has two important implications for interpretation.
First, when demographic growth is exponential and effectively unbounded, as in
the geometric model with constant growth rate \(n_t \equiv n\), the
non-compounding budget equation reproduces the familiar peak-and-collapse
behaviour of classical Ponzi schemes: high promised coupons combined with rapid
expansion create a transient accumulation of capital, followed by inevitable
depletion once inflows slow or stop. In this sense, classical Ponzi schemes can
be captured as special parameter constellations within the present framework.

Second, once the demographic law is modified to feature a finite pool of
potential investors and a decreasing effective growth rate---as in the
quasi-logistic and SIR-based specifications---the same budget equation admits
parameter regimes in which liabilities shrink in line with participation while
the asset side continues to earn market returns. In those regimes, the system
behaves as a finite-horizon pooled-income investment product rather than as a
genuine Ponzi scheme: all investors receive the promised coupons and the
repayment of \(I_0\), and the promoter ends the cycle with a non-negative
capital stock at termination. The subsequent analysis explores these regimes in
detail.


\subsection{Traffic light scenarios: an intuitive approach}\label{Sect3p2}

\begin{figure*}[p!]  
\vspace{-0.0cm}\centering
\includegraphics[width=0.75\linewidth]{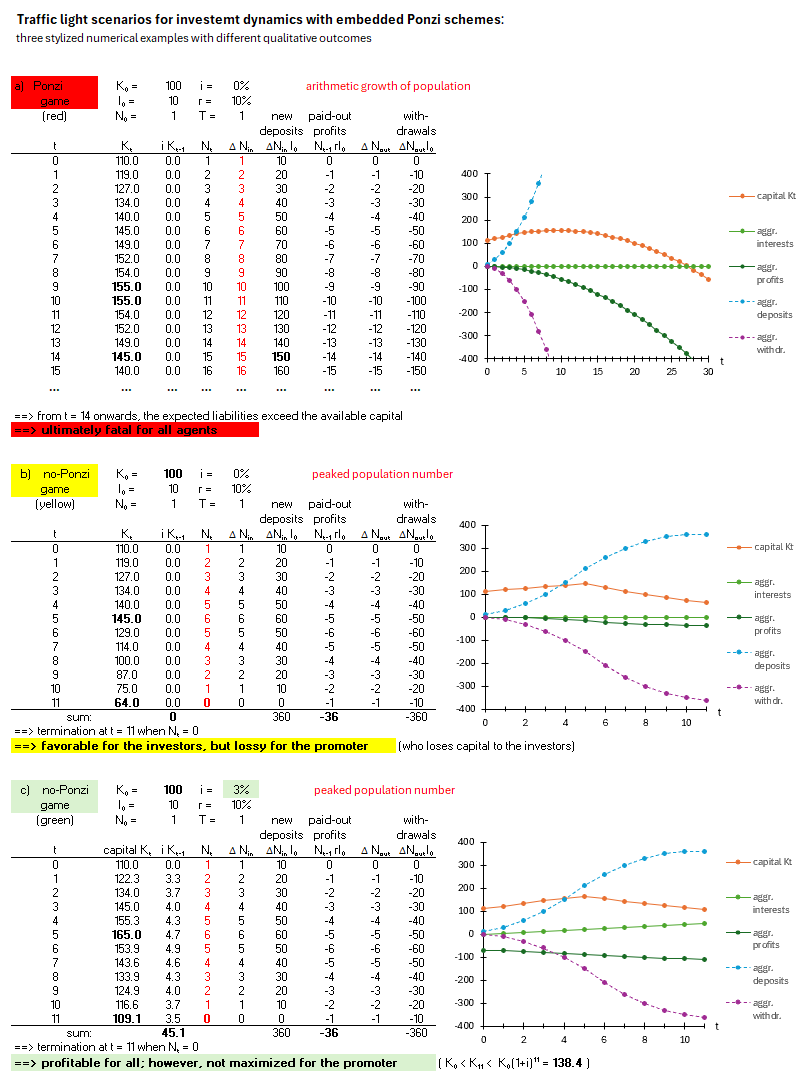}
\caption{Traffic light scenarios: The evolution of capital and its associated contributions within a rudimentary investment framework is examined. The numerical dynamics can be traced sequentially. The interaction among the initial capital, compounded interest, the number of investors, and the particular lock-up period \(T = 1\) yields qualitatively distinct outcomes. The accompanying graphs depict the temporal evolution of the capital and its aggregated contributions. The red light scenario corresponds to a final collapse (top panel), the yellow light scenario encounters profiting investors but a promoter with some losses (middle panel), and the green light scenario only sees agents with gains at termination (bottom panel). }
\label{Fig3_TrafficLightScenarios}
\end{figure*}

Before introducing sophisticated models that require more advanced mathematics, this section adopts an intuitive step-by-step approach, laying the groundwork for some of the fundamental results presented later. As in any investment context, there is no free lunch: payouts must ultimately be financed by genuine revenues. These revenues typically derive from interest earned on the initial capital as well as from compounded investment returns. Importantly, the relevant interest rates may be significantly lower than the overall investment rates of return ($i < r$). The promised profits therefore translate into relatively high dividend distributions — a feature frequently observed in high-yield stocks in financial markets. However, such dynamics cannot be sustained if the number of investors continues to increase indefinitely. Consequently, \emph{the combination of relatively short lock-up periods and a demographic profile that peaks over time can yield alternating phases of asset growth and decline, with a positive net balance in the end}.

The contrived examples shown in Figure \ref{Fig3_TrafficLightScenarios} provide transparent traceability of cash flows within such a system. The capital evolution obeys the budget equation \eqref{Ktrec}. From panel to panel, all but one of the parameter values remain the same. The lock-up time is set to one period ($T=1$) in all panels, meaning that each participant remains invested for only one period before leaving the system permanently (hence $N_t=\Delta N_{in}$).

In part a) (top panel), a classic Ponzi scheme emerges, characterized by an arithmetically increasing number of participants and an initial phase of steady capital growth. Around \(t=9\), however, the system reaches its peak. At this point, either the promoter absconds with the remaining funds, or the available cash begins a sharp and continuous decline. Outstanding obligations can no lon\-ger be met, inevitably driving the system into bankruptcy (around \(t=14\)). Ultimately, all investors — including the promoter, who loses the initial capital — suffer losses. This scenario is marked with a red traffic light.  

In part b) (middle panel), the number of investors follows a humped trajectory: after reaching its peak (at about $t=5$), participation gradually declines until all investors have withdrawn their deposits and exited the system. During this process, the promoter’s initial capital is partially absorbed by the investors. Participants receive both a return on their investment and, eventually, the repayment of their initial deposit. For a limited period, investing in the system proves more profitable than investing in the capital market, since \(r > i\). However, the promoter made a  loss. This scenario is marked with a yellow traffic light.  

In part c) (bottom panel), once the peak has been reached, both the number of investors and the outstanding obligations decline in parallel, supported by interest generated on the aggregated capital. There is no paradox here; the system functions as a competitive balance between aggregated inflowing interest and aggregated outflowing profits. Investors behave rationally — seeking to maximize profits — and with integrity, acting in a moral\-ly impeccable manner. Participation occurs within a transitory scheme managed by a promoter who is not a con artist but rather an altruistic funder, content with earning submaximal profits on his own behalf. This scenario is marked with a green traffic light.

In summary, the three numerical examples illustrate qualitatively different possible outcomes of capital dynamics (including Ponzi-like behavior), referred to as \emph{traffic light scenarios}: 
\begin{enumerate}
\item[-] The red scenario corresponds to an ultimately collapsing system (\(K_t < 0\) at termination time \(t\)). 
\item[-]  The yellow scenario reflects a harmed promoter but successful investors ($0\le K_t < K_0^\mathrm{pro}$ and \(I_0 \cdot (1 + rT) > I_0$). 
\item[-] The green scenario represents a benign outcome for all participants ($0<K_0^\mathrm{pro}<K_t<K_0^\mathrm{pro}\cdot(1 + i)^t$ and $I_0 \cdot (1 + rT) > I_0$).
\end{enumerate}


\vspace{0.2cm}
\subsection{Capital dynamics with geometric growth (constant rate $n$)}\label{Sect3p3}
\vspace{0.4cm}
We begin our examination of capital dynamics in specific investment systems by employing the geometric population dynamics introduced in Section \ref{Sect2p1} and embedding them into the capitalist framework specified by equation \eqref{Ktrec}. We first consider a simplified scenario characterized by persistent participation (i.e., investments without withdrawals, corresponding to the limit $T \to \infty$), arriving at the recurrence relation
\begin{eqnarray}
\Delta K_{t} & =& i K_{t - 1} + \Delta N_{t}^{in}\ I_{0} - N_{t - 1}rI_0\\
K_{t} &=& K_{t - 1} +\Delta K_{t} \nonumber\\
&=&(1 + i)K_{t - 1}+ I_{0}N_{0}(n - r)\, (1 + n)^{t - 1} 
\label{GMKtimplicit}
\end{eqnarray}

\subsubsection{Unlimited investment with persistent participation}\label{Sect3p3p1}

We distinguish between two types of persistent participation, those that are not and those that are exposed to the capital market.

(i) \emph{No access to the capital market ($T \rightarrow \infty$, $i = 0$).} Setting for now $K_0^{pro}=0$ and allowing neither for market returns ($i=0$) nor for withdrawals ($T\rightarrow \infty$), the new investors fully finance the profits of the other onboarded investors with their investments: payments (deposits) of the new investors equals payouts (profits) to the existing investors, hence $\Delta K_t=0$ and therefore $\Delta N_{t}^{in}I_{0}= N_{t - 1}P$ or $N_{0}{(1 + n)}^{t - 1}n\ I_{0} = N_{0}{(1 + n)}^{t - 1}rI_{0}$; this yields the simple \emph{no-Ponzi-game (NPG) condition}
\begin{equation}
n \ge r \, . 
\label{NPG}
\end{equation}
As long as the growth rate \(n\) of system participants is at least as large as the percentage return \(r\) promised to investors, the system functions either as a zero‑sum game in the case \(n = r\), or it generates a temporary surplus through capital accumulation when \(n > r\). However, in practice the influx of new participants inevitably stagnates at some point. Once \(n < r\), the deposits of new members are no longer sufficient to sustain payouts at the established levels, rendering the system loss‑making and ultimately leading to collapse. Consequently, the so‑called NPG condition can never be satisfied permanently. Referring in this situation to a NPG solution is therefore inaccurate and misleading. Stop and loss, or collapse, and with it significant harm to numerous investors, is ultimately unavoidable.

\begin{figure}[!t]
\centering
\includegraphics[width=.99\linewidth]{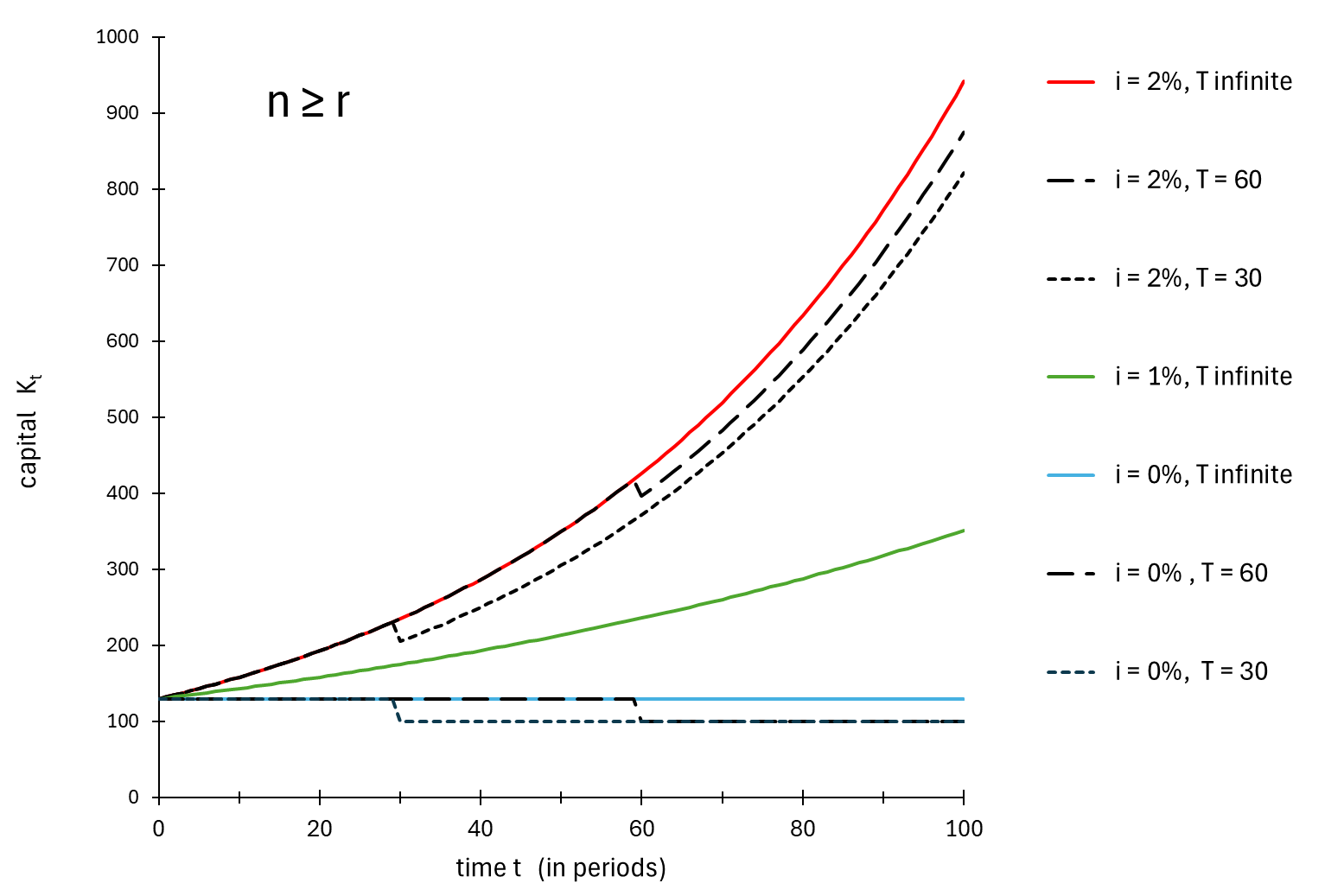}
\includegraphics[width=.99\linewidth]{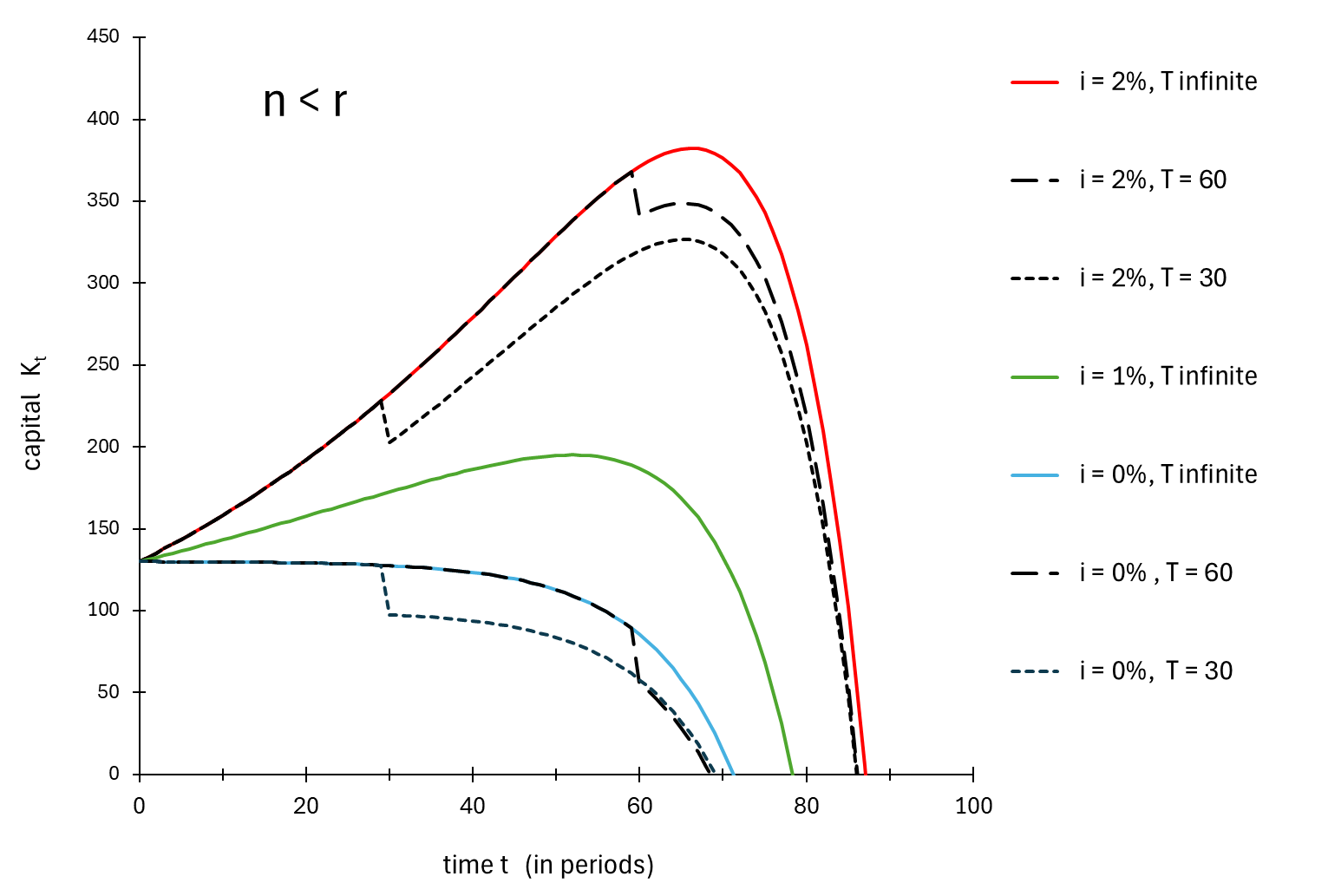}
\caption{Capital dynamics with geometric growth ($n$ = constant). Parameter values used: constant rate of return on investment $r = 10\%$ and growth rate $n = r$ in the upper panel, but $n = 9.95\% <r$ in the lower panel. Market-based rate of return $i$ and duration $T$ of investments as indicated in the legends. Initial values $K_{0}^{pro}= 100$, $I_{0} = 3$, $N_{0}= 10$.}
\label{Fig4}
\end{figure}

In Figure \ref{Fig4} a typical evolution of an asset without access to the capital market $(i = 0\%$, $T$ infinite) is depicted with blue solid lines, once for $n\ge r$ (upper panel) and once for $n<r$ (lower panel).

(ii) \emph{With exposure to the capital market ($T \rightarrow \infty$, $i \geq 0$)}. From equation \eqref{GMKtimplicit}, one infers the NPG condition: the available capital will always grow if $n>r$; if $n<r$, capital may initially increase due to interest on the initial endowment, but once $\Delta K_t<0$ and hence $ iK_{t - 1}<I_{0}N_{0}(r-n)\, (1 + n)^{t - 1}$, it will decline and ultimately drop to zero. By applying the prolific theorem presented in Section \ref{Sect3p3p3} to the implicit relation \eqref{GMKtimplicit}, one readily arrives at the explicit solution:
\begin{equation}
K_{t} = \left\{ K_{0} - N_{0}I_{0}\frac{n - r}{n - i} \right\}{(1 + i)}^{t} + {\ N}_{0}I_{0}\frac{n - r}{n - i}{(1 + n)}^{t}.
\label{GMKttsmallerT}
\end{equation}
(This holds for any time $t<\infty$; the paramter $T$ is irrelevant here.) For the case $n$=$i$=$\,0$ the capital becomes $K_t=K_0-N_0I_0r\,t$ \footnote{Taking the $\lim_{n,i\rightarrow 0}((1+n)^t-(1+i)^t)/(n-i)\approx \lim_{n,i\rightarrow 0}((1+tn)-(1+ti))/(n-i)=t$.} with collapse time $t_{collapse}=K_0/(N_0I_0r)$. A Ponzi scheme that assumes geometric growth in order to describe the evolution of the number of incoming and redeemed
investors results in exponential cash inflow and outflow, which may seem unrealistic \cite{Peng21}. It never\-theless captures the baseline idea of competing capital flows and delivers a mathematically tractable, albeit rather simple, approach to the phenomenon of a concave graph for the capital evolution. This may explain the motivation of an unscrupulous operator to set up such a scheme in order to intentionally disappear with the accumulated net capital when it peaks. Formulas for the calculation of both the peak time and the time of collapse will be provided in the next subsection.

\subsubsection{Limited duration of an investment due to a lock-up time}\label{Sect3p3p2}

The lock-up period $T$ is an investor's waiting time until withdrawal of the investment. Different to the previous subsection with continuous participation, here each investor remains invested for only a finite number of periods and then exits the system. This occurs for the very first time at $t = T$. The impact on the capital dynamics formally presents as follows. As long as $t<T$ the capital stock evolves according to equation \eqref{GMKttsmallerT}. At $t=T$ first repayments by the first exiting investors with a total amount of $N_0I_0$ are due, thus
\begin{equation}
	K_T = \left\{ K_0 - N_0I_0\frac{n - r}{n - i} \right\}{(1 + i)}^{T} + {\ N}_{0}I_{0}\frac{n - r}{n - i}{(1 + n)}^{T} - N_0I_0.\label{GMKT}
\end{equation}
Then the explicit solution for all periods is\footnote{\emph{Proof}: For $t>T$, the capital at the beginning of a period is given by the capital of the period before (including market-based returns) $+$ the deposits of new investments $-$ the payouts to existing investors $-$ the repayments to exiting investors. Formally, by means of the corresponding implicit solution $K_{t} = K_{t - 1} + \Delta K_{t} = (1 + i)K_{t - 1} +\Delta N_t^{in}\ I_0- N_{t - 1}\, P -\Delta N_t^{out}I_0$ ---with $N_t$, $\Delta N_t^{in}$, and $\Delta N_t^{out}$ as given in equations \eqref{GNt}, \eqref{GNtin}, and \eqref{GNtout}, respectively--- one has
\begin{equation*}
K_{t} = (1 + i)K_{t - 1} +N_{0}\ I_{0}\ (n - r)\left\lbrack 1 - (1 + n)^{-T} \right\rbrack{\cdot (1 + n)}^{t - 1}.
\end{equation*}
The corresponding explicit solution readily follows from the prolific theorem given in equations \eqref{KtRR1}f., with $t_0=T$ and $K_T$ from \eqref{GMKT}. $\square$}
\begin{equation}
\boxed{%
K_{t} =
	\begin{cases} 
	\left( K_{0} - {\ N}_{0}I_{0}^{}\frac{n - r}{n - i} \right){(1 + i)}^{t} + {\ N}_{0}I_{0}\frac{n - r}{n - i}{(1 + n)}^{t}
	& t<T\\
	 \left( K_{T} - {\ N}_{0}I_{0}\frac{n - r}{n - i}\eta (1 + n)^{T} \right)(1 + i)^{t - T} &\\ 
	\hspace{3.3cm}+N_{0}I_{0}\frac{n - r}{n - i}\eta{(1 + n)}^{t} & t \ge T
\end{cases}
}
\end{equation}
(with the constant factor $\eta\equiv 1 - (1 + n)^{-T}$). The times for peak capitalization and for the completed collapse are found by means of the conditions $\Delta K_{t}=0$ and $K_t=0$, respectively, delivering
\begin{eqnarray}
\hspace{-0.5cm}
t_{collapse}&=&\frac{1}{\log\left( \frac{1+n}{1+i}\right)}\times \nonumber\\
& &\hspace{-1.5cm}\begin{cases}
\log\left( 1+\frac{n-i}{r-n}\cdot\frac{K_0}{N_0I_0}\right) & (T\rightarrow \infty) \\
\log\left( 1+\frac{n-i}{r-n}\cdot\frac{K_T}{N_0I_0}\cdot \frac{1}{(1+n)^T-1}\right)+T& (T<\infty),
\end{cases}\\
 &=& t_{peak} + \Delta t_{pop}.
\label{tcoll}
\end{eqnarray}
Herein, the duration of the pre-disaster phase, i.e., the timespan from peak until collapse ---called here \emph{period of precipice}\footnote{Borrowing from a narrative or dramaturgical context, this disastrous final stretch might also be called a pre-catastrophic stage, eve of disaster, or the moments of last suspence.}--- is found to be for both cases 
\begin{equation}
\Delta t_{pop} = \log\left( \frac{n}{i}\right) / \log\left( \frac{1+n}{1+i}\right) \,\,\,-1.
\end{equation}
For a peak and a collapse to occur, the following \emph{peak-and-crash condition} must be fullfilled:
\begin{equation}
r>n>i>0.
\end{equation}
If $i=0$, one cannot calculate $\Delta t_{pop}$ because the capital does not peak but instead decreases monotonically from the initial value to zero. As a numerical example, in Figure \ref{Fig4} (lower panel with $n<r$), the uppermost curve peaks at 66.7 periods, enters a precipice phase lasting 24.1 periods, and reaches the time axis after 86.1 periods, based on the parameter set $\{K_0, I_0, N_0, i, n, r\} = \{130, 10, 3, 0.02, 0.0995, 0.1\}$.

\subsubsection{On a prolific theorem concerning particular recurrence equations}\label{Sect3p3p3}

First-order linear difference (or recurrence) equations with constant coefficients and exponential inhomogeneities admit a simple solution for discrete times \(t \in \mathbb{N}\). This solution was mentioned prominently several times in the previous text.

\noindent
\fbox{%
  \parbox{\columnwidth}{%
\emph{\textbf{Theorem}}: The recurrence equation 
\begin{equation}
K_{t} = (1 + i)K_{t - 1}\ \  + c\,(1 + n)^{t - 1}, \hspace{0.5cm}(t > t_{0} \geq 0)\label{KtRR1}
\end{equation} 
with initial value $K_{t_0}$ at $t=t_0$ and a real parameter $c$, admits the explicit solution
\begin{equation}
K_{t} = \left( K_{t_0} + \frac{c}{i - n}{(1 + n)}^{t_{0}} \right)(1 + i)^{t - t_{0}} - \frac{c}{i - n}(1 + n)^{t}.\label{Kt1}
\end{equation}
  }%
}

\vspace{0.3cm}This useful result can be verified either by substituting the explicit expression \eqref{Kt1} into the implicit relation \eqref{KtRR1} or by performing a complete induction. Alternatively, following a standard approach, the general solution can be derived via the method of characteristic polynomials. It is then obtained as the sum of the homogeneous solution and a particular solution, subject to the given initial conditions, as worked out in the Notes.\footnote{\emph{Proof} of the prolific theorem (equations \eqref{KtRR1}-\eqref{Kt1}): The textbook treatment for solving an inhomogeneous linear difference equation of the form $x_n=a\,x_{n-1}+c\,b^{n-1}$, with initial value $x_{n_0}$ (for $n>n_0$), consists in determining its homogeneous solution $x_n^{(h)}$ and additionally finding some particular solution $x_n^{(p)}$, with their sum composing the general solution. For the homogeneous solution (setting $c=0$), the Ansatz $x_n=r^n$ provides the characteristic polynom $r^n-ar^{n-1}=0$, with solution $r=a$, hence $x_n^{(h)}=Ca^n$, where the value of the constant $C$ will follow from some initial conditions. Similarly, inserting the Ansatz $x_n^{(p)}=D\,b^n$ for the particular solution into the original difference equation, the characteristic polynom leads to $D=c/(b-a)$, yielding $x_n^{(p)}=c/(b-a)\cdot b^n$ (requiring $b\neq a$).  Hence the general solution is $x_n=x_n^{(h)}+x_n^{(p)}=Ca^n+c/(b-a)\cdot b^n$. Starting with $n=n_0$ as the initial condition, the constant becomes $C=a^{-n_0}\left(  x_{n_0}-c/(b-a)\cdot b^{n_0}\right)$. Thus, $x_n=\left( x_{n_0}-c/(b-a)\cdot b^n \right)a^{n-n_0}+c/(b-a)\cdot b^n$. Identifying $t=n$, $t_0=n_0$, $K_t=x_n$, $a=1+i$, and $b=1+n$, equation \eqref{Kt1} is recovered. The case $b=a$ corresponding to $n=i$ is omitted here; for computational purposes, simply set $n=i+\epsilon$, with $\epsilon\ll 1$. \,\,$\square$}

\emph{Special case} \(t_{0} = 0\): The recurrence relation
\begin{equation}
K_{t} = (1 + i)K_{t - 1}\ \  + c\, (1 + n)^{t - 1}, \hspace{0.5cm}(t > 0)\label{KtRR2}
\end{equation}
with initial value \(K_{0}\), has the solution
\begin{equation}
K_{t} = \left( K_{0} + \frac{c}{i - n} \right){(1 + i)}^{t} - \frac{c}{i - n}{(1 + n)}^{t}.\label{Kt2}
\end{equation}

\emph{Corollary}: If in equation \eqref{KtRR1} there is not only a single inhomogenei\-ty of the form $c (1 + n)^{t - 1}$  but there are multiple inhomogeneities given as a sum $\sum_{j = 1}^{k}c_{j}\left( 1 + n_{j} \right)^{t - 1}$, then relations \eqref{KtRR1} to \eqref{Kt2} hold subject to obvious corresponding adjustments. 

Additional applications of this theorem and its related corollary concerning classical Ponzi schemes are discussed in \cite{Par13} and \cite{Par17}. Moreover, the following triple of published models, which either feature iterative or sum-based solutions, can also be explicitly solved using our formalism; for instructional purposes, these solutions are derived in the Notes. (i) The model in \cite{Ati21} allows for the realization of nominal interests with rate $r_n$ and for withdrawals in terms of a fraction of the currently available deposits. Their solution for the capital dynamics depends on a sum that, equivalently, can be replaced with our formalism by an explicit formula.\footnote{In particular, the entries for the cash inflow and cash outflow in the budget equation for the capital,  $K_t=(1+r_n)K_{t-1}+D_{t-1}-W_{t-1}$, are progressively growing deposits $D_t=(1+r_i)^t$ (with growth rate $r_i$ that acts similar to a demographic growth rate) and withdrawals $W_t$ that are a fraction $r_w$ of the deposits and become a finite geometric series. Concretely, withdrawals are found to be a difference of progressions, $W_t=\beta(1+r_i)^{t-1}\sum_{k=0}^{t-1}[(1+\alpha)/(1+r_i)]^k  =(W_0+\beta/(\alpha-r_i))(1+\alpha)^t-(\beta/(\alpha-r_i))(1+r_i)^t$, with $1+\alpha = (1+r_p)(1-r_w)$ and $\beta=r_w(1+r_p)D_0$. The solution for $K_t$ relies on their theorem 2.1 for a first-order nabla difference equation (related to Sturm-Liouville theory) and involves a sum with $t$ terms, namely, $K_t= K_0(1+r_n)^t+(1+r_n)^{t-1} \sum_{s=0}^{t-1}(1+r_n)^s (D_s-W_s)$ (with $W_0=0$, their eq. 3.4). Instead, adopting (the corollary of) our prolific theorem one straightforwardly and explicitely finds \[K_t =(K_0-T_1+T_2)(1+r_n)^t+T_1(1+\alpha)^t-T_2(1+r_i)^t,\] with $T_1=[W_0+\beta/(\alpha-r_i)] / (r_n-\alpha)$ and $T_2=[D_0+\beta/(\alpha-r_i)] / (r_n-r_i)$.} (ii) The "simple, stylistic Ponzi model" suggested in \cite{BP25} (their Section 4.2.2) in order to form some intuition on fraudulent schemes takes into account constant periodic takeouts by the promoter.\footnote{The formal treatment is as follows: the cash balance for the capital $K_t$ at discrete times $t$ reads $K_t=K_{t-1}+I_t-O_t-C_t$, where the cash inflow and outflow are assumed to be $I_t =I_{t-1}(1+g)=I_0(1+g)^t$ (with growth rate $g$) and $O_t=I_{t-1}(1+r)$ (with promised rate of return $r$), respectively, and starting with $K_0=I_0$. The periodic takeout by the promoter is $C_t=c$ with $C_0=0$. Then --by means of the prolific theorem-- the implicit relation $K_t=(1+i)K_{t-1}+I_0(g-r)(1+g)^{t-1}-c(1+h)^{t-1}$ (with ad hoc rates $i$ and $h$ that later are to become zero) has the explicit solution $K_t=\{K_0-I_0(g-r)/(i-g)+c/(i-h)\}(1+i)^t-I_0(g-r)/(i-g)\cdot(1+g)^t-c/(i-h)\cdot(1+h)^t$.  The term proportional to $c$ is  $[(1+i)^t-(1+h)^t]/(i-h)$; taking $h=i+\epsilon$, it reduces to $t$ after taking the required limits $\lim_{i\rightarrow 0}$ and $\lim_{\epsilon\rightarrow 0}$. Hence the solution becomes \[ K_t=I_0(g-r)/g\cdot(1+g)^t-ct+I_0r/g,\] with indeed $K_0=I_0$ at $t=0$. This is a superposition of an increasing exponential sequence (assuming $g>r$) and of a decreasing arithmetic sequence. See \cite{BP25} for a discussion.} (iii) A spreadsheet model for a classical Ponzi scheme is presented in \cite{Win21}, with the capital evolution being subject to periodic partial withdrawals that are proportional to the aggregated liabilities (including promised ones and hence some fictitious asset). Once again, the prolific theorem enables a direct and analytic approach to the problem.\footnote{In \cite{Win21}, the capital fund at time $t$ (that is, at the end of period $t$, but \emph{before} the addition of new deposits and the distribution of payouts to investors, termed "start cash") evolves from its value in the previous period ($K_{t-1}$), decreased by the payouts made ($P_{t-1}$), and increased by new deposits ($N_0I_0$). The sum of these is then compounded over one period at the risk-free interest rate $i$. This provides the recurrence relation $K_t = (1+i)(K_{t-1}+I_0N_0-P_{t-1})$, starting with $K_0$. Herein, payouts $P_t=wL_t$ are considered to be partial withdrawals and calculated as some amount proportional to the overall liabilites $L_t$, with $w$ the fraction withdrawn. Liabilities include both investments and promised return on investments, with a promised rate of return $r$. In the beginning of each period, payouts reduce the liabilities. The evolution of the liabilities thus obeys the recurrence equation $L_t=(L_{t-1}-P_{t-1})(1+r)+N_0I_0=(1-w)(1+r)L_{t-1}+N_0I_0$, with initial value $L_0=N_0I_0$. Applying the prolific theorem, one finds $L_t=N_0I_0/(\lambda-1)\,(\lambda^{t+1}-1)$, with $\lambda\equiv(1-w)(1+r)$. Inserted into the recurrence equation for the capital, the (corrolary to the) prolific theorem results in \[K_t=[K_0+c_1/(i-\lambda+1)+c_2/i](1+i)^t-[c_1/(i-\lambda+1)]\lambda^t-c_2/i,\] with $c_1\equiv -[(1+i)w/(\lambda-1)]N_0I_0$ and $c_2\equiv -[r(1-w)/w]c_1$.}


\subsection{Capital dynamics with decreasing growth rates $n_t$}\label{Sect3p4}

In Section \ref{Sect3p2}, it became clear that the ratio $n/r$ is essential for the capital dynamics of the Ponzi scheme described there: for $n/r \ge 1$, sufficient capital remains in the system or even continues to grow. In the latter case, however, this only lasts until the pool of potential participants is exhausted. For $n/r < 1$, the capital itself always becomes exhausted after a certain time. The present section analytically describes mixed scenarios in that they may begin with a favorable ratio $n_t/r \ge 1$, but due to a decreasing growth rate $n_t$ eventually exhibit a fatal ratio $n_t/r < 1$. As discussed below, both the quasi-logistic growth model and the non-standard SIR model exhibit this behaviour.\footnote{The logistic capital management theory or financial satu\-ration theory (e.g., \cite{Kar25}) can reproduce the same qualitative behavior when its Ansatz is suitably extended. In its canonical form, the model assumes that capital $K_t$ grows proportionally to a logistically increasing number of investors $N_t$, i.e., $K_t=I_0N_t=K_\text{sat}/[1+A\,(1+j)^{-1}]$, where $j$ is the growth rate and $A=(K_\text{sat}-K_0)/K_0$. Thus, capital evolves from its initial value $K_0$ toward the saturation level $K_\text{sat}$. To capture a realistic decline phase, one may instead adopt the bi-logistic Ansatz of the Jolicoeur-Pontier-model, consisting of a growing and a declining logistic term in the denominator, \[K_t=K_\text{sat}\,[A\,(1+j)^{-t}+B\,(1+k)^{t}]^{-1},\] where $k$ denotes the rate of decline (typically $B$$\ll$$A$ and $k$$>$$j$). This extension naturally produces a final downturn after some growth phase, capturing the effect of investors withdrawing their capital.}


\subsubsection{Quasi-logistic growth}\label{Sect3p4p1}

We again distinguish the cases $t< T$ before any withdrawals (that is equivalent to setting $T\rightarrow \infty$) and $t \ge T$ with withdrawals occuring. 

(i) $t< T$. If the budget equation $K_{t} = K_{t - 1} + \Delta K_t=(1 + i)K_{t - 1} +\Delta N_t^{in} I_0 - N_{t - 1}^{t < T}\ rI_0 $ for the system's capital is based on the quasi-logistic population dynamics described in Section \ref{Sect2p2}, one obtains the recurrence relation
\begin{eqnarray}
K_{t} &=& (1 + i)K_{t - 1} + I_{0}N_{t - 1}^{t < T}\ (n_{t} - r)\label{QLDeltaKta}\\
&=& (1 + i)K_{t - 1} + I_{0}N_{0}\frac{1}{n\left( 1 - \frac{N_{0}}{N} \right)} \,\times \nonumber \\
& &\hspace{1.5cm}{(1 + n)}^{t - 1}n_{t - 1}(n_t - r)\label{QLDeltaKtb}
\end{eqnarray}
Herein, the decreasing growth rate $n_t$ as defined in equation \eqref{ntQL} has an initial value close to $n$. Obviously, as long as $n_t>r$ the aggregated capital increases but will decrease when $n_t<r$. Starting with $K_0$ at $t=0$, after a few formal iterations the pattern becomes apparent, and a compact notation for the explicit solution $K_t ^{QL}\equiv K_t$ is found to be
\begin{equation}
\boxed{K_t ^{QL}= (1 + i)^{t}\left(  K_{0} + I_{0}N_{0} \,S_t^{QL} \right),
\hspace{0.5cm}(1\le t < T)}
\label{QLKttsT}
\end{equation}
with the decisive sum formula
\begin{equation}
S_t^{QL} = \sigma\sum_{k = 1}^{t }\left( \frac{1 + n}{1 + i} \right)^{k-1} \,n_{k-1} (n_{k} - r)\label{QLSt}
\end{equation}
and the ad hoc constant $\sigma = \left[\left( 1 - \frac{N_{0}}{N} \right)n(1 + i) \right]^{-1}$; by definition, $S_0^{QL} \equiv0$.

In the limit $N \rightarrow \infty$ and hence with constant values $n_{k-1} = n_{k} = n$, the well-known sum formula for finite geometric series allows to straightforwardly transform this solution into the one for the geometric model (equation \ref{GMKttsmallerT}). In case of abstinence from the market (i.e., $i=0\%$), the hump-time condition $\Delta K_t=0$ becomes equivalent to $n_t-r=0$, yielding in turn the hump time, i.e., the time of maximal capitalization:
\begin{equation}
t_{K_t\,peak}^{QL}\vert _{_{i=0\%}}=\frac{\log \left[ \left( \frac{N}{N_0}-1\right) \left( \frac{n}{r}-1\right) \right]}{\log (1+n)}.
\label{KttpeakQL}
\end{equation}
For $i>0\%$, this characteristic period can instead be found numerically and always occurs some time later as in the market unfriendly case with $i=0\%$.

(ii) $t\ge T$. If one requires every investor to leave the system (and to withdraw the deposit) after the expiration of an investment time-span of $T$ periods (a.k.a. lock-up period), the system's downsizeing process becomes accelerated, as expected and as formally described as follows. At $t=T$ the solution given in \eqref{QLKttsT} needs to be corrected by the repayments to the very first and now leaving investors, numbered $\Delta N_T^{out}=\Delta N_{T-T}^{in}=N_0$ according to equation \eqref{QLDNtout}, hence 
\begin{equation}
K_T^{QL}= (1 + i)^{T}\left\{  K_{0} + I_{0}N_{0} \, S_T^{QL}  \right\} -I_0N_0 .\label{QLKT}
\end{equation}
Starting at $t=T$ with an initial value $K_T^{QL}$, the budget equation for the system's capital afterwards is similar to equation \eqref{QLDeltaKta} with an additional withdrawal term, i.e.,  $K_t= (1 + i)K_{t - 1} + \Delta N_{t}^{in}\ I_0 - N_{t - 1}\ rI_0- \Delta N_{t - T}^{in}\ I_0$ (be\-ing equivalent to equ. \ref{Ktrec}) , and reads  
\begin{eqnarray}\hspace{-0.5cm}
K_{t} &\hspace{-0.4cm}=&\hspace{-0.4cm} (1 + i)K_{t - 1} +I_0 \left\{  N_{t - 1}^{t<T}(n_t-r) -N_{t -T-1}^{t<T}(n_{t-T-1}-r)\right\}\nonumber \\
        &\hspace{-0.4cm}=&\hspace{-0.4cm} (1 + i)K_{t - 1} +I_{0}N_{0}\frac{1}{n\left( 1 - \frac{N_{0}}{N} \right)}{(1 + n)}^{t - 1}\times\nonumber\\
        &  & \hspace{-0.3cm} \Big\{ n_{t - 1}(n_t - r) - (1+n)^{-T}n_{t -T- 1}(n_{t-T} - r)\Big\} .
\end{eqnarray}
After a few formal iterations, the explicit solution $K_t ^{QL}\equiv K_t$ for times $t\ge T$ is inferred to be
\begin{eqnarray}
K_t^{QL}&=&(1+i)^{t-T}\Bigg\{  K_T^{QL} \nonumber \\
& &\hspace{-0.8cm}+\,\,I_0N_0 \Big[  (1+i)^T\left( S_t^{QL}-S_T^{QL} \right) - S_{t-T}^{QL} \Big] \Bigg\},
\label{KtQL}
\end{eqnarray}
Inserting $K_T^{QL}$ (equ. \eqref{QLKT}) and $S_x^{QL}$ (equ. \eqref{QLSt}) yields

\begin{empheq}[box=\fbox]{align}
K_t^{QL}&=(1+i)^{t}\Bigg\{  K_0 +\,\,I_0N_0 \Big[   S_t^{QL}\Big.\nonumber \\
& \hspace{-0.8cm}- (1+i)^{-T}  \left( 1+S_{t-T}^{QL}\right)  \Big] \Bigg\}.\hspace{0.5cm}(t\ge T)
\label{QLKttgeT}
\end{empheq}
\newline

\begin{figure}[!t]
\centering
\includegraphics[width=.99\linewidth]{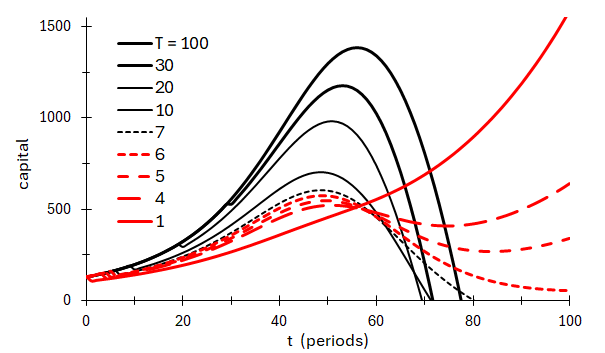}
\caption{Quasi-logistic investment dynamics: capital evolution for a Ponzi scheme-like model with a decreasing demographic growth rate $n_t$. Initial values  are $K_{0}^{pro}= 100$, $I_{0} = 3$, $N_{0}= 10$, and fixed rates $n=10\%$, $r=5.2\%$, and $i=3\%$. Examples for schemes with quasi-logistic growth of member numbers (using formulas \eqref{QLKttsT} and \eqref{QLKttgeT}). The investment timespan (or lock-up period) $T$ varies as indicated, with the case $T=100$ periods representing here continued participation. Ultimately collapsing systems are shown in black ($T\ge 7$), while red lines represent surviving systems ($T\le 6$).}
\label{Fig5}
\end{figure}

\begin{figure}[t!]
\centering
\includegraphics[width=.88\linewidth]{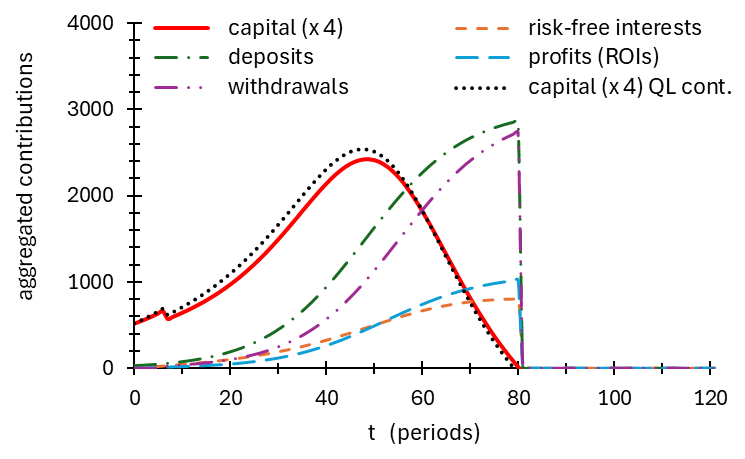}
\includegraphics[width=.88\linewidth]{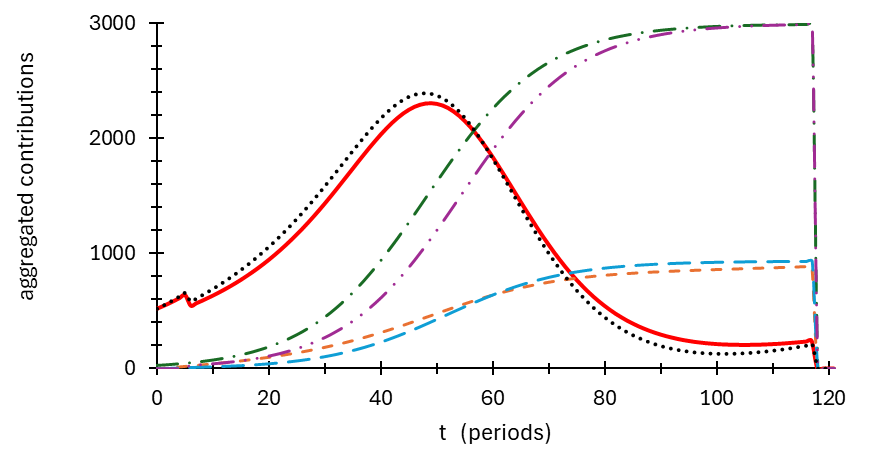}
\includegraphics[width=.88\linewidth]{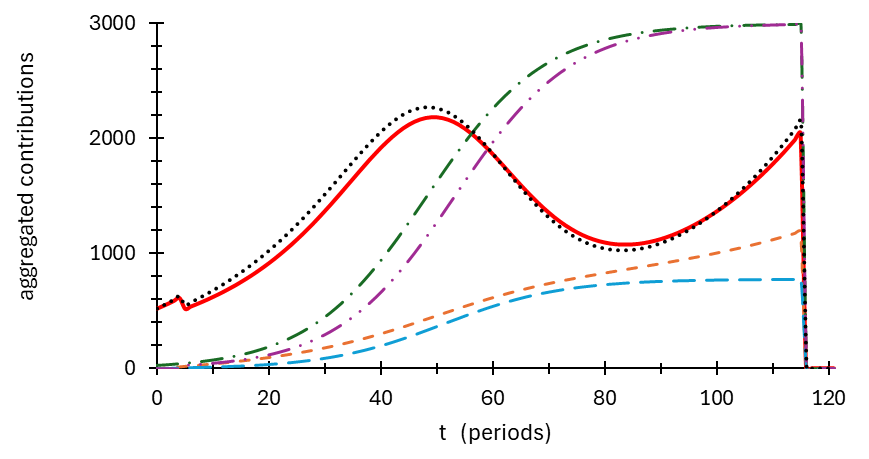}
\includegraphics[width=.9\linewidth]{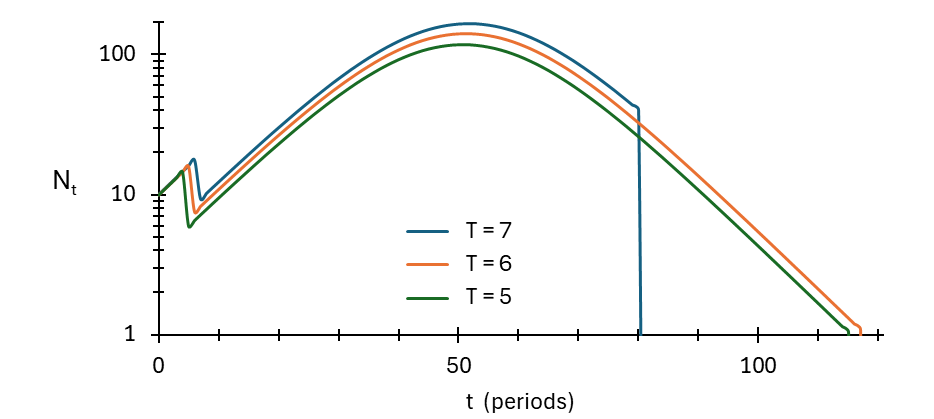} 
\caption{Critical capital dynamics for the quasi-logistic investment scheme. From Fig.~\ref{Fig5}, the cases for $T=7,6,5$ are extracted and shown again in the respective top three panels. The available capital $K_t$ (scaled by a factor of 4 for visualization purposes) is shown by solid red lines. It is supplemented by aggregated quantities: risk-free interest earned on capital $\sum iK_t$, inflows from deposits of new investors $\sum\Delta N_t^{\rm in}I_0$, outflows due to pro\-fits $\sum N_{t-1}rI_0$ (i.e., promised returns on investment, ROI), and withdrawals $\sum\Delta N_t^{\rm out}I_0$. For $T=7$, the system collapses at period 81 ($K_t<0$). In contrast, for  $T=6$ and $T=5$, the schemes do not collaps but terminate calmly when the number of participants drops below 1. The final capital positions of the promoter are: $K_t<K_0^{pro}$ (i.e., not profitable) for $T=6$ and $K_t>K_0^{pro}$ (i.e., profitable) for $T=5$. These three distinct outcomes correspond to different traffic light scenarios (cf. Section \ref{Sect3p2}). The black dotted lines trace the capital for the related continuous-time model presented in the Appendix (with adjusted values $r =$ 6.1\%, 6.2\%, 6.3\% ceteris paribus for the cases $T=7,6,5$, respectively). The bottom panel shows the evolution of the number of investors $N_t$ over time.}
\label{Fig6}
\end{figure}

Capital evolution under quasi-logistic demographics (and thus a bounded pool of investors) with a sufficiently short lock-up period \(T\) behaves as a non-Ponzi, finite-horizon pooled investment game rather than as a classical Ponzi scheme. For finite \(N\) and decreasing \(n_t\), the quasi-logistic model generates SIR-like humps in participation and, under suitable parameter choices \((i,r,T,\ldots)\), admits benign no-Ponzi paths in which all investors are repaid and the promoter ends with a positive capital surplus at termination. By contrast, for \(n_t \equiv n\) and \(N \to \infty\), the model collapses to the geometric Ponzi-scheme framework with exponential participation and the familiar peak-and-crash dynamics; a similar outcome may arise for unfavourably long lock-up periods. These phenomena are illustrated in Figure~\ref{Fig5}. In that figure, only the lock-up period \(T\) is varied, as indicated in the legends. Typically, the system eventually collapses (black lines). However, for sufficiently short lock-up periods (red lines), a characteristic feature of models with declining population growth emerges: given some modest initial capital \(K_0^{\mathrm{pro}}\) provided by the promoter and a modest positive risk-free market rate \(i>0\) (even in cases with \(r>i\)), the system need not collapse for suitably small maturities \(T\). In such cases, market returns outweigh coupon payments to an ever-decreasing number of participants. Once the last investor exits and no one has incurred a loss, the promoter receives the accumulated remaining capital as compensation for successfully operating the scheme under legal conditions. Restricting attention to the pair \((i,T)\) and treating these as variables, any scheme that numerically satisfies the no-Ponzi condition
\begin{equation}
K_t^{QL}(K_0^{\mathrm{pro}},i,T)\big\vert_{I_0,N_0,N,n,r} > 0
\label{KiiT}
\end{equation}
for all times \(t > T\) operates in a no-Ponzi environment. The multivariate function in~\eqref{KiiT} thereby defines a positive-valued, higher-dimensional no-Ponzi surface. A general analytical characterization of such no-Ponzi conditions, however, lies beyond the scope of the present paper.\newline

To further illustrate the no-Ponzi game (NPG) or survival effect, the cases shown in Figure \ref{Fig5} for investment timespans of $T=$ 7, 6, and 5 periods are highlighted again in Figure \ref{Fig6}, accompanied by relevant aggregated quantities. The bottom panel displays the number of current investors $N_t$, which exhibits the typical hump characteristic of growth models with effectively decreasing growth rates. Evidently, the longer the investment timespan (lock-up period), the greater the number of investors present and the higher the peak of this hump.

The top three panels depict the available capital $K_t$ as red solid lines. Additionally, aggregated quantities are presented, including interests on capital ($\sum iK_t$), inflow from deposits of newly joining investors ($\sum \Delta N_t^{in}I_0$), and outflow ($\sum \Delta N_t^{out}I_0$), among others. Due to decreasing growth rates $n_t$, the number of members reaches a peak before monotonically declining. This behavior is reflected in the saturation of both aggregated inflowing cash from deposits and aggregated outflowing cash from withdrawals. 

In the top panel (for $T=7$), obligations can no longer be satisfied after about $t=81$ periods, leading the system to collapse. In contrast, in the second and third panels (for $T=6$ and $T=5$, respectively), the declining investor base critically reduces payments of returns on investments (ROIs). These are offset by income derived from effective risk-free interest in the capital market. Over time, the market-based yields eventually exceed the aggregate payments to participants, whose number steadily diminishes. Once the final investor exits the system, the scheme terminates naturally, leaving the promoter with a residual asset — either a loss (second panel) or a profit (third panel).  In the latter case, however, the residual value is suboptimal: it is lower than the profit that would have accrued had the promoter simply invested the initial capital $K_t^{\mathrm{pro}}$ in the capital market under compound interest. Thus, although all agents — promoter and investors alike — exit the scheme with profits, the promoter falls short of the maximum attainable profit, whereas all investors realize their maximum possible gains. In the terminology of the traffic-light classification introduced in Sect.~\ref{Sect3p2}, this outcome corresponds to the green-light scenario.

Taken together, the formalism and examples presented so far support the existence of collapse-free scenarios that yield benign outcomes and, in some cases, profits for all participating agents. The aim of this paper is not to engage with the complex theoretical challenge of formulating a multivariate no-Ponzi game criterion involving up to eight variables, as the graphical results already convey sufficient insight. 

\subsubsection{Nonstandard SIR-model ($T_0=0$)}\label{Sect3p4p2}

The population dynamics of a non-standard SIR mo\-del ---recently proposed in \cite{Lem25}, equipped with explicit formulas, and recalled here in Section \ref{Sect2p3p2}--- serves as the basis for the development of the system's capital at all times $t\ge 0$. In particular, the respective number of investors in the system ($N_t$, being identical with $I_t$ in the original notation but renamed here to fit with our notation), those entering the system ($\Delta N_t^{in}$), and those quitting the system ($\Delta N_t^{out}$) are
\begin{eqnarray}
 N_t&=&N_0\left( \frac{1+\beta}{1+\gamma} \right)^t p_t,\\ 
\Delta N_t^{in}&=&\beta\frac{N_{t-1}}{N_{t-1}+S_{t-1}}S_t \nonumber\\
&=&\beta S_0 \left( 1+\frac{S_0}{N_0} \left( \frac{1+\gamma}{1+\beta} \right)^{t-1}\right)^{-1}p_t,\\
\Delta N_t^{out}&=&\gamma N_t= \gamma N_0\left( \frac{1+\beta}{1+\gamma} \right)^t p_t.
\end{eqnarray}
Therefore, the net change of capital at time $t$ due to inflowing interest income and incoming new investments and due to outflowing returns and withdrawals is formally expressed by $\Delta K_t = i K_{t-1}+\Delta N_t^{in} I_0 - N_{t - 1}rI_0 -\Delta N_t^{out} I_0$ and becomes
\begin{equation}
\Delta K_t =i K_{t-1}+I_0N_0\, b_t\,p_t,
\end{equation}
with time-dependent coefficient
\begin{equation}
\hspace{-0.3cm}
b_t= \left( 
\frac{\frac{S_0}{N_0} -r\left( \frac{1+\beta}{1+\gamma} \right)^{t-1}}{\frac{S_0}{N_0} +\left( \frac{1+\beta}{1+\gamma} \right)^{t-1}}\,\beta -r-\gamma\frac{1+\beta}{1+\gamma} \right) \left( \frac{1+\beta}{1+\gamma} \right)^{t-1}
\label{bt}
\end{equation}
and the product function $p_t$ as given by equation \eqref{pt}. Setting $\gamma >0$, as in \cite{Lem25}, withdrawals are assumed to occur at \emph{all} times after the initial time. (This corresponds to setting $T_0=0$ in an extended model that will be investigated in the next subsection.) Starting at $t=0$ with $K_0$ and performing a few subsequent formal iterations $K_t=K_{t-1}+\Delta K_t$, one finds the generating pattern and hence the general result
\begin{equation}\label{KtnsSIR}
\boxed{K_t^{nsSIR}=(1+i)^t\left( K_0+I_0N_0\,S_t^{nsSIR}\right)},
\end{equation}
where the second term contains a decisive summation over the product of $b_k$ and $p_k$,
\begin{equation}\label{StnsSIR}
S_t^{nsSIR}=\sum_{k=1}^t \frac{b_k\,p_k}{(1+i)^k}.
\end{equation}

\begin{figure}[!t]
\centering
\includegraphics[width=.82\linewidth]{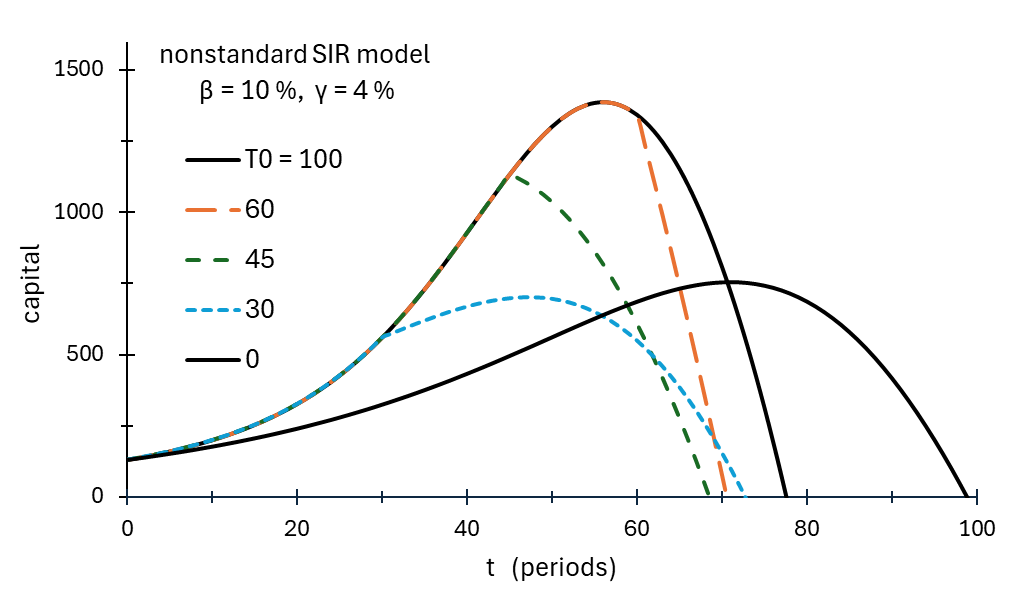}
\includegraphics[width=.82\linewidth]{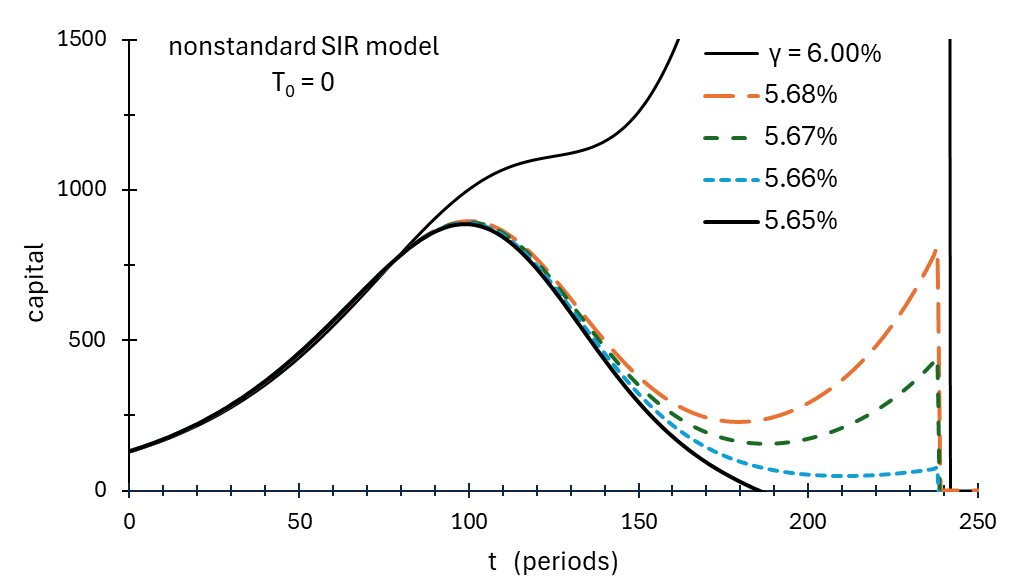}
\includegraphics[width=.82\linewidth]{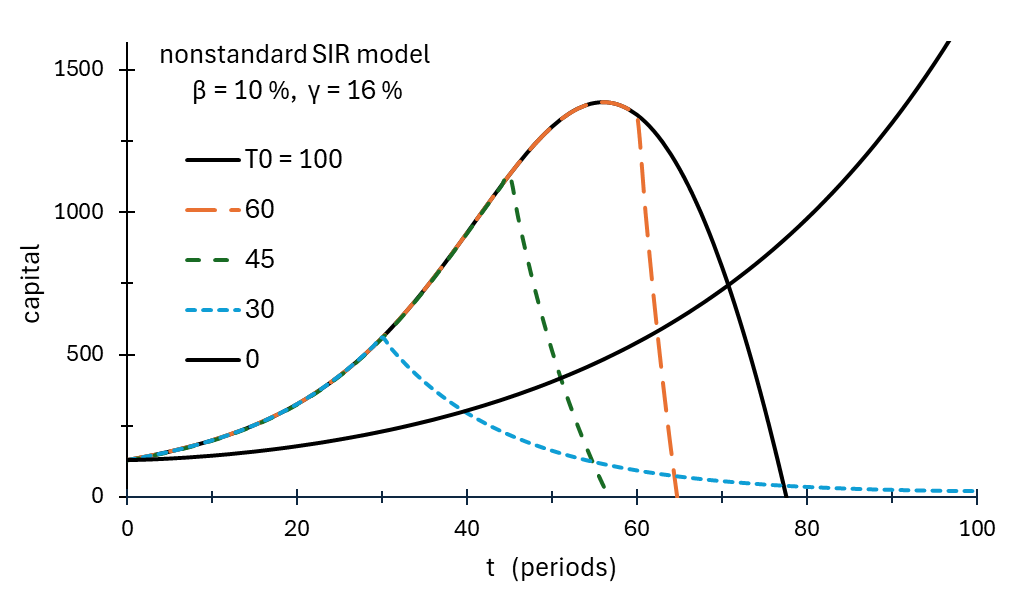}
\includegraphics[width=.82\linewidth]{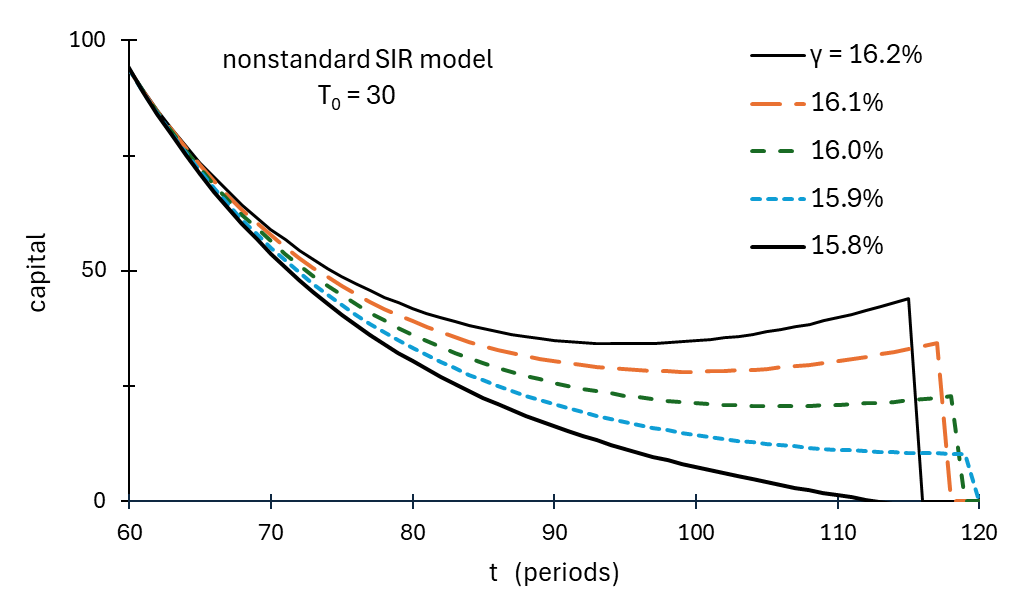}
\caption{Capital evolution in an investment system governed by nonstandard SIR demography with delayed initial recovery. As in Figures \ref{Fig4} and \ref{Fig5}, initial values  are $K_{0}^{pro}= 100$, $I_{0} = 3$, $N_{0}= 10$, with fixed rates $r=5.2\%$, $i=3\%$. The examples shown are based on formulas \eqref{KtnsSIR2} and \eqref{KtnsSIRT0}, assuming a joining rate $\beta = 10\%$ (equal to $n$) and a quitting rate $\gamma =0\%$ until $t=T_0$, and $\gamma > 0$ afterwards, as indicated. In each panel, only one parameter is varied, as indicated in the legends. The case $T_0=0$ corresponds to the ordinary nonstandard SIR model (i.e., with immediate onset of recovery), while $T_0=100$ represents continued participation (i.e., indefinitely delayed recovery). In the top panel, all systems represent red-light scenarios with Ponzi collapse. The scenario with $T_0=0$ in the top panel is shown again in the panel below, with slightly different values for $\gamma$. The scenario with $T_0=30$, as depicted in the third panel from the top, is contained again in the bottom panel, together with affiliated cases. Vertical drop of lines occurs when the number of participants has declined to below one, i.e., when the system's operation time is over.}
\label{Fig7}
\end{figure}

\subsubsection{Nonstandard SIR-model with $T_0\ge 0$}\label{Sect3p4p3}

If initial recoveries or withdrawals by quitting investors do not occur from the very start of the epidemic or system operation but only begin after a delay of $T_0>0$ periods, then the formalism presented in the previous subsection splits into two parts. This division arises from the demographic description in Section \ref{Sect2p3p3} and corresponds to the following straightforward adaptations.

(i) $1\le t\le T_0$ (withdrawals are omitted yet, hence the recovery rate is $\gamma =0$): Starting with $K_0$ and $N_0$ (and $R_0=0$ and hence $S_0=N-N_0$), the solution found previously (equations \eqref{KtnsSIR}-\eqref{StnsSIR}, together with \eqref{bt} for $b_k$ and \eqref{pt1} for $p_k$) readily delivers
\begin{equation}\label{KtnsSIR2}
\boxed{K_t^{nsSIR}=(1+i)^t\left( K_0+I_0N_0\,S_t^{nsSIR,\,<}\right)},
\end{equation}
where the decisive summation over the product of $b_k^{<}$ and $p_k^{<}$,
\begin{equation}
S_t^{nsSIR,\,<}=\sum_{k=1}^t \frac{ b_k^{<}\,p_k^{<} }{(1+i)^k},
\end{equation}
goes with the adapted time-dependent coefficient
\begin{equation}
b_k^{<} = \left( 
\frac{\frac{S_0}{N_0} -r\left( 1+\beta \right)^{k-1}}{\frac{S_0}{N_0} +\left( 1+\beta \right)^{k-1}}\,\beta -r
 \right)
\left( 1+\beta \right)^{k-1}
\label{bttsT0}
\end{equation}
and with the product function $p_k^{<}$ as given by equation \eqref{pt1}. No withdrawals occur within this first interval, these only start with period $T_0+1$ at the beginning of the second interval. 

In the special case of perpetual participation with no withdrawals (i.e., $T_0 \rightarrow \infty$ or at least $t<T_0$ with $T_0$ sufficiently large) and when the capital is not invested in an external market ($i=0\%$), the peak condition $\Delta K_t=K_t^{nsSIR}-K_{t-1}^{nsSIR}=0$ implies $b_t^{<}p_t^{<}=0$, indicating that the capital peak occurs at
\begin{equation}
t_{K_t\ peak}^{nsSIR,<}\vert _{i=0\%}=\log \left( \frac{\beta-r}{r}\,\frac{S_0}{I_0}\right)\,/\,\log (1+\beta),
\label{tpeaknsSIRT01}
\end{equation}
given $\beta >r$. For $\beta <r$, there's no peak but a monotonic decline to collapse at $K_t=0$. If $i>0\%$, the peak must be looked for numerically or graphically.

(ii) $t\ge T_0+1$ (with regular withdrawals enforced, corresponding to taking $\gamma >0$): Starting the iterations with initial values $K_{T_0}$ and $N_{T_0}$ and giving a similar formal treatment as before, the general result for times $t>T_0$ becomes 
\begin{equation}
\boxed{K_t^{nsSIR}=(1+i)^{t-T_0}\left( K_{T_0}+I_0N_{T_0}\,S_t^{nsSIR,\,>}\right)},
\label{KtnsSIRT0}
\end{equation}
where now the second term contains the decisive summation 
\begin{equation}
S_t^{nsSIR,\,>}=\sum_{k=T_0+1}^t \frac{b_k^{>}\,p_k^{>}}{(1+i)^{k-T_0}}.
\end{equation}
Herein enter the time-dependent coefficient
\begin{equation}
b_k^{>}=\left(
\frac{\frac{S_{T_0}}{N_{T_0}} -r\left( \frac{1+\beta}{1+\gamma} \right)^{k-T_0-1}}{\frac{S_{T_0}}{N_{T_0}} +\left( \frac{1+\beta}{1+\gamma} \right)^{k-T_0-1}}\,\beta
-r-\gamma\frac{1+\beta}{1+\gamma} 
\right) \left( \frac{1+\beta}{1+\gamma} \right)^{k-T_0-1}
\label{bttgT0}
\end{equation}
and the product function $p_k^{>}$ as given by equation \eqref{pt2}. 

In the special case of delayed initial withdrawals (i.e., with $T_0 < \infty$) and with the capital not invested in an external capital market ($i=0\%$), the peak condition $\Delta K_t=K_t^{nsSIR}-K_{t-1}^{nsSIR}=0$ implies $b_t^{>}p_t^{>}=0$ or $b_t^{>}=0$, indicating that the capital peak now occurs at
\begin{equation}
t_{K_t\ peak}^{nsSIR,>}\vert _{i=0\%} = T_0+
\frac{ \log\left(  
                   \frac{ \beta-r(1+\gamma)-\gamma}{(r(1+\gamma)+\gamma)(1+\gamma)  }\frac{S_{T_0}}{I_{T_0}}
	\right)
       } 
       { \log \left( \frac{1+\beta}{1+\gamma} \right) }.
\label{tpeaknsSIRT02}
\end{equation}
This holds if $\beta>r(1+\gamma)-\gamma$ and if anyway the argument of the nominator's $\log$-term is larger than one. Otherwise, if the argument equals one, the peak occurs at $T_0$, and else there's no peak at all but a monotonic decline from the inital $K_0$ to collapse (i.e., to $K_t=0$). Again, if $i>0\%$ the peak has to be found numerically or graphically. Equation \eqref{tpeaknsSIRT01} is recovered for $\gamma=0$ (implying $T_0=0$, too).

Setting $T_0=0$, formula \eqref{KtnsSIR2} (or \ref{KtnsSIRT0}) describes a Ponzi system with (unmodified) nonstandard SIR-model demography, as dealt with in the previous subsection \ref{Sect3p4p2}. This significant case is included in the upper three panels of Figure \ref{Fig7} as one of the black solid lines: if the capital hits the time axis ---that is, when it reaches zero---, the system crashes ($\beta>\gamma\ge 0$). But, remarkably, the capital also may take off ($\beta<\gamma$). Similar scenarios can be seen for $T_0\ge0$ in general. These different outcomes show ---as analogously already encoutered with capital dynamics based on the quasi-logis\-tic growth model--- that a crash may not only be obstructed, but avoided at all under certain predetermined circumstances. Depending on the parameter values, the capital will either be depleted before all investors exit the system, causing the Ponzi scheme to collapse, or — after a finite operating period when all investors have exited — the promoter will be left with some remaining capital, resulting in a harmless and productive conclusion to the no-Ponzi investment scheme. A relevant, implicit premise is the exposure of the capital to the external market ($i>0$); setting $i=0\%$ inevitably leads to bankruptcy. 

The Ponzi models in \cite{Peng21} and \cite{BP25} rest on extended SIR-models that, in addition to our approach with the nonstandard SIR-model, incorporates lock-up periods, too. Numerically solving their coupled delayed differential equations should similarly exhibit, for adequately chosen parameter values, the sustainability effect. 

\subsubsection{A look at a recurrent NPG investment scheme}\label{Sect3p4p4}

The issue of sustainability is also briefly addressed in two other studies cited above in the literature review of Section 1.

(i) "A simple (but unusual) condition under which the Ponzi scheme may never go bank\-rupt" is stated in \cite{Parlar25}. Translated into our terminology, the condition given in their Theorem 6.1 (for a continuous model, however) recalls the fact that if the exponential demographic growth rate \(n\) (and hence the contribution or deposit rate) exceeds the promised rate of return \(r\), i.e., \(n > r\), then the cash position will grow monotonically. A similar statement can be formulated for their discrete-time model\footnote{Indeed, assuming a geometric sequence $u_k=u_0(1+\hat{r})^k$ for the net deposits by the investors (with $\hat{r}$ being the growth rate) and inserting it according to their Theorem 2.1 into their equation (3) yields for the cash position of the Ponzi fund after $t$ periods the explicit solution \[c_t=c_0+u_0\{ (1-r/\hat{r})[(1+\hat{r})^t-1]-(rs_0/u_0-r/\hat{r})\,t \}.\] This is ---analogous to Note \#16--- a superposition of a geometric sequence and an arithmetic sequence. Therefore, for $n\equiv\hat{r}>r$ the term exponential in $t$ will dominate over the one linear in $t$.}.  This is equivalent to our no-Ponzi game condition \eqref{NPG} and hence relies on the unrealistic assumption of an infinite pool of participants.

(ii) In their game theoretic model with two overlapping generations, the authors in \cite{Bar25} (their Appendix B) consider "a variation of the model in the text in which the imposter [i.e., the fraudulent promoter, called long-lived agent] can undertake a risky but profitable investment." They go on with "In contrast to our benchmark model in which the long-lived agent only benefits from stealing and will necessarily default on some cohort [i.e., short-lived investors], adding investment will make it possible for the long-lived agent to avoid default." However, while their construction of a 'Ponzi equilibrium' parallels our findings of a non-defaulting system, in their model, riskiness is associated with an undisclosed high interest rate ($i>r$ in our notation), in contrast to our model, where even the case ($i<r$) can be sustainable. 

We thus propose a novel path to sustainability, derived from our dynamic investment model. In a system where all investors have left one has $N_t=0$. Notwithstanding, in the frameworks with decreasing demographic growth rates the real number of participants at late times is a real null sequence. Modelling within finite times thus requires to stop the calculation at some predetermined value $N^\star \equiv N_{t^\star} \gtrsim 0$ corresponding to some end-of-efficacy or termination time $t^\star$. Picking out the QL-formalism and solving equation \eqref{QLNt} algebraically for $t$  (at times $t\ge T$), one finds
\begin{equation}
t^\star (N^\star) = \frac{1}{1+n} \cdot
\begin{cases}
\log u^{-}   & t \le t_{N\,peak} \\
\log u^{+}  & t \ge t_{N\,peak},
\end{cases} 
\end{equation}
where the peak time $t_{N\,peak}$ stems from equation \eqref{QLtNpeak} and where
\begin{eqnarray}
u^{\pm} &=& \frac{-B\pm\sqrt{B^2-4AC}}{2A},
\end{eqnarray}
with
\small
\begin{eqnarray*}\hspace{-0.5cm}
A &=& \frac{N^\star N_0}{N^2(1+n)^T}\\
B &=&-\left( 1-\frac{N_0}{N} \right) \left( 1-\frac{1}{(1+n)^T} - \frac{N^\star}{N} \left(1+\frac{1}{(1+n)^T}\right)  \right) \\
C &=&\frac{N^\star}{N_0}\left( 1-\frac{N_0}{N}\right)^2 .
\end{eqnarray*}
\normalsize
\begin{figure*}[p!]
\centering
\includegraphics[width=0.85\textwidth]{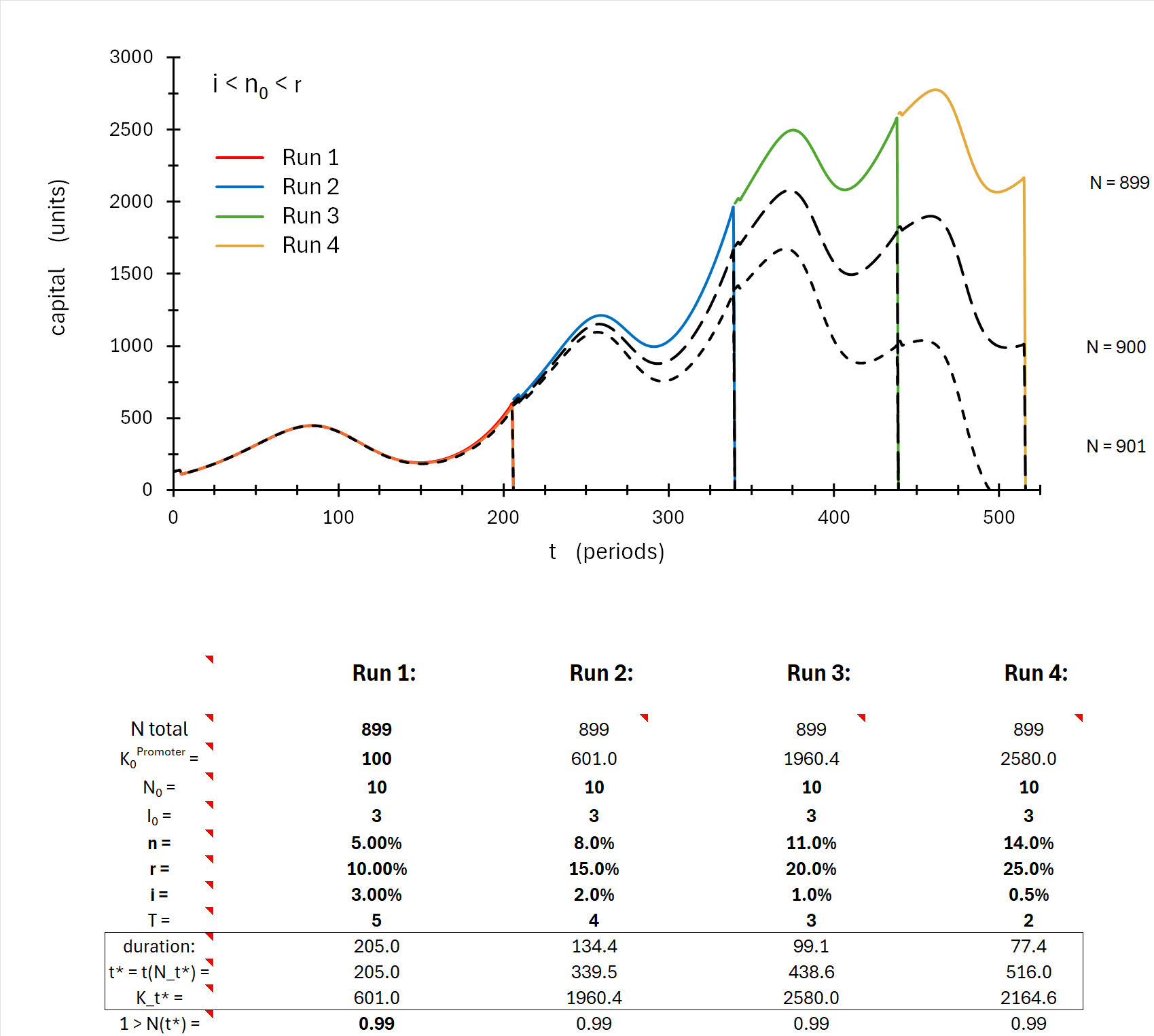}
\caption{The feasibility of recurrent profitable runs in a suitably constructed investment system is illustrated by the solid lines. Specific parameter values for each of the four runs are provided in the table below the graph. The system’s high sensitivity to most parameters is demonstrated by the dashed curves, where the initial size of the participant pool is held constant across simulations. This pronounced dependence on parameter configurations makes the operator’s task particularly challenging, as values must be carefully recalibrated for each run.}
\label{Fig8_recurrent}
\end{figure*}
A sustained no-Ponzi game concludes with enforced termination at a stopping time $t^\star$, when, say, $N^\star=0.99$ (i.e., when there is less than one investor left). At this point, the promoter holds capital $K_{t^\star}>0$, which may serve as the initial endowment $K_0^{\text{pro}}$ for a subsequent instance of the investment game. The  vector of  parameter values either remains unchanged as compared to the previous run or is adjusted according to the conditions for a next game. For practical purposes, each run evolves on the interval beginning at the (relative) time $0$ and terminating at a (relative) stopping time $t^\star$. Consequently, the analytical framework introduced for the QL-model can be applied repeatedly across successive runs. Subsequently, each run is positioned at the appropriate location on the global or absolute time axis, yielding a piecewise construction of the overall trajectory. Figure \ref{Fig8_recurrent} illustrates such a sequence of instances, with the varying sets of parameter values given in the accompanying table. Actually, only  a few of the eight main parameters experience some changes from run to run, namely, $n$, $r$, $i$, and $T$, the others are kept fix or adjusted. For the first three runs, the graphs exhibit an uptrending pattern of expansion and contraction. In the last, fourth run, the system's capital records a large loss, mostly due to insufficient market activity. Furthermore, the dashed lines stress the extreme sensitivity of the model's behaviour due to changes in the total number of potential investors, $N$ (or of the ratio $N_0/N$).

In a practical financial setting, an operating team must retain control over key variables and carefully orchestrate the evolution of capital flows. After each run, model parameters can be recalibrated to reflect shifting demographic or market conditions, as well as newly defined objectives. The option to pause operations altogether also remains part of the strategic toolkit.

The main difference from traditional banking lies in the emphasis on a lock-up period that is short enough to ensure system stability. However, the application of maturity periods is by no means alien to banking practice, as evidenced by derivative instruments such as convertible securities and other options. By requiring a controlled number of investors as another key element, our model may therefore serve as a foundation for constructing an innovative derivative that mirrors the capital evolution presented here.

Madoff’s house of cards endured for an unusually long time. In light of our dynamic model, which encapsulates Ponzi‑like systems, the existence of a recurrent investment game with considerable longevity appears plausible —  and thus potentially attractive —  though it remains operationally very challenging. Nevertheless, improper triggering or unpredictable exogenous shocks may still bring the house down —  that is, cause the system to collapse —  thereby resembling the downfall of a classic Ponzi scheme, with no prospect of a bailout.


\section{Conclusions}\label{Sect4}

This paper investigates the capital dynamics of three mathematically related investment models — geometric, quasi-logistic, and SIR-based — each admitting\linebreak closed-form solutions. Within this framework, classical Ponzi-type behaviour appears as a special case arising from adverse parameter constellations, while favourable configurations yield sustainable and profitable investment processes. The quasi-logistic formulation, in particular, synthesizes essential features of the geometric and SIR-based approaches.

A central determinant of system behaviour is the dynamic ratio $n_t/r$, which separates regimes of capital growth and decline. Both constant and time-varying cases were analysed, including scenarios where $n_t/r \ge 1$ at first but later falls below~1, as seen in the quasi-logistic and non-standard SIR models. Features such as finite investment horizons, delayed withdrawals, or optimized quitting rates show that even markets with moderate interest rates may allow sustainable operation.

A \emph{profitable but transitory} investment mechanism requires initial capital, market exposure, a declining investor growth rate, and short lock-up or delayed recovery times. Such systems typically exhibit a hump-shaped evolution of investor numbers and end with positive net capital — often exceeding the initial capital —  assuming a promoter satisfied with reduced profit, such as a mission-driven or philantropic seed investor.

Three qualitative “traffic-light” outcomes emerge:\linebreak \emph{red} for ultimate collapse (final capital $K_t<0$), \emph{yellow} when investors profit but the promoter incurs losses ($K_0^{pro}>K_t\ge 0$), and \emph{green} when outcomes are beneficial for all stakeholders (though suboptimal for the promoter). No-Ponzi game (NPG) patterns display a rise and fall of both participation and capital while remaining non-negative, enabling repeated cycles with increasing initial capital. Unlike “rational” infinite-agent roll\-over constructions \cite{OCZ92}, the systems discussed here are finite-time, parameter-sensitive investment games that can represent either deceptive Ponzi structures or transparent, legally admissible mechanisms, including PAYG-type analogues.\footnote{\label{fn20}The model does not incorporate regular deposits over $T$ or overlapping generations.}

Taxes and fees, though excluded from the core model, would reduce the effective interest rate $i$ and return rate $r$, potentially reversing inequalities such as $i<n<r$. Numerical tests show that many viable NPG scenarios persist under such modifications, illustrating the limited relevance of simple analytical conditions like~\eqref{NPG} and the inherent uncertainty in the relationship between $r$ and $n$.\footnote{\label{fn21}Cf.\ the heuristic argument in \cite{Man15} regarding $r-\tau<n$.}

Given the eight underlying parameters, no general analytical NPG condition is available; numerical investigation is required. All explicit formulas were verified against iterative solutions, and spreadsheet implementations used for figures are available on request. Because the recursive method recalculates all quantities period by period, it can naturally accommodate parameter changes or stochastic variation \cite{Peng21,Parlar25,Par14}. A conti\-nuous-time version of the quasi-logistic model (i.e., taking the limit of infinitesimaly short periods) is nevertheless provided in \ref{AppB}.

A key qualitative result is the emergence of recurrent NPG cycles with fluctuations around an upward trend. These investment cycles — characterized by transient expansions followed by controlled contractions  —  closely resemble \emph{mean-reversion patterns} and \emph{asset price bubbles} observed in financial markets. This connection places the present framework within the wider literature on boom–bust dynamics and recurrent bubbles \cite{OCZ92,Dom15,Per21,Bar25,Hir24}, where transient, self-limiting expansions are seen as fundamental economic phenomena.

Possible applications include short-term high-yield funds, certain governmental bond structures, and social-security systems, the latter requiring extensions to overlapping generations. Finally, the models highlight the smooth continuum between legitimate banking operations and fraudulent Ponzi-like systems: small parameter changes can shift otherwise lawful structures into illegitimate territory. As emphasized in \cite{Basu14}, camouflaged Ponzi mechanisms can be difficult to detect because they often intertwine with legitimate activities; in practice, such investment tools should therefore be accompanied by adequate securities.


$\newline$
\small

\paragraph{\bf Acknowledgments} This paper was originally motivated by an intended---but ultimately unrealized---contribution to the Spring 2025 philosophy workshop on conflicts between values and interests, offered by Dr.\ Suzann-Viola Renninger (University of Zurich, Switzerland), with a particular emphasis on Ponzi schemes. The author gratefully acknowledges the public and free accessibility of both the ETH Library (Swiss Federal Institute of Technology, Zurich) and the Zurich Central Library (Zentralbibliothek). The author also recognizes arXiv—originally created by Paul Ginsparg and now operated by Cornell University with major support from the Simons Foundation and a global network of member institutions—as an indispensable tool for the literature review. Perplexity.ai assisted in refining the language of the manuscript.

\normalsize


\theendnotes







\normalsize


\appendix

\section{Formal iterations} 
For the demographic evolution under geometric and quasi-logistic growth, explicit iterations of the recurrence relations are carried out to derive the corresponding analytical formulae.

\newpage 
\subsection{Demography with geometric growth}
\begin{table}[h!]\label{AppendixTable1}
\centering
\tiny
\begin{tabular}{|c|p{25mm}|c|c|}  
\hline
$\mathbf{t}$ & $\mathbf{N_t}$ & 
$\mathrm{\Delta}\mathbf{N_t^{in}}$ & \
$\mathbf{\mathrm{\Delta}}\mathbf{N_t^{out}}$ \\
\hline
0 & \(N_{0}\) & \(N_{0}\) & 0 \\
1 & \(N_{0}(1 + n)\) & \(N_{0}n\) & 0 \\
2 & \(N_{0}{(1 + n)}^{2}\) & \(N_{0}(1 + n)n\) & 0 \\
\vdots & \vdots & \vdots& \vdots \\
\(\mathbf{t}\) &
\(\mathbf{N}_{\mathbf{0}}\mathbf{(1 + n)}^{\mathbf{t}}\) &
\(\mathbf{N}_{\mathbf{0}}\mathbf{(1 + n)}^{\mathbf{t - 1}}\mathbf{n}\) &
\textbf{0} \\
\vdots & \vdots & \vdots& \vdots \\
\(T - 1\) & \(N_{0}{(1 + n)}^{T - 1}\) & \(N_{0}{(1 + n)}^{T - 2}n\) & 0 \\
\hline

\(T\) & \(N_{0}(1 + n)^{T} - N_{0}\) \newline
\({= N}_{0}(1 + n)^{T}\left\lbrack 1 - \frac{1}{{(1 + n)}^{T}} \right\rbrack\)
& \(N_{0}{(1 + n)}^{T - 1}n\) & \(N_{0}\) \\

\(T + 1\) & \(N_{0}{(1 + n)}^{T} - N_{0}\) \newline
\({+ N}_{0}(1 + n)^{T}n - N_{0}n\) \newline
\({= N}_{0}(1 + n)^{T + 1} - N_{0}(1 + n)\) \newline
\({= N}_{0}(1 + n)^{T + 1}\left\lbrack 1 - \frac{1}{{(1 + n)}^{T}} \right\rbrack\)
& \(N_{0}{(1 + n)}^{T}n\) & \(N_{0}n\) \\

\(T + 2\) & \({= N}_{0}(1 + n)^{T + 1} - N_{0}(1 + n)\) \newline
\({+ N}_{0}{(1 + n)}^{T + 1}n - N_{0}(1 + n)n\) \newline
\({= N}_{0}(1 + n)^{T + 2} - N_{0}{(1 + n)}^{2}\) \newline
\(= N_{0}(1 + n)^{T + 2}\left\lbrack 1 - \frac{1}{{(1 + n)}^{T}} \right\rbrack\) \newline
& \(N_{0}{(1 + n)}^{T + 1}n\) & \(N_{0}(1 + n)n\) \\
\vdots & \vdots & \vdots& \vdots \\
\(\mathbf{t}\) &
\(\mathbf{N}_{\mathbf{0}}\mathbf{(1 + n)}^{\mathbf{t}}\left\lbrack \mathbf{1 -}\frac{\mathbf{1}}{\mathbf{(1 + n)}^{\mathbf{T}}} \right\rbrack\)
& \(\mathbf{N}_{\mathbf{0}}\mathbf{(1 + n)}^{\mathbf{t - 1}}\mathbf{n}\)
&
\(\mathbf{N}_{\mathbf{0}}\mathbf{(1 + n)}^{\mathbf{t - T - 1}}\mathbf{n}\) \\
\hline
\end{tabular}
\normalsize
\caption{Evolution of the number of members in the geometric growth model (with constant growth rate $n$). The number of all active investors $N_{t}$, the number of new(ly entering) investors $\Delta N_{t}^{in}$, and the number of exiting investors $\Delta N_{t}^{out}$ evolve iteratively according to equations \eqref{recrelNt}, \eqref{DNin}, and \eqref{GDNout}, respectively.}
\end{table}

\subsection{Demography with quasi-logistic growth}
\begin{table}[h!]\label{AppendixTable2}
\centering
\tiny
\begin{tabular}{|c|p{25mm}|p{23mm}|c|}  
\hline
$\mathbf{t}$ & $\mathbf{N_t}$ & $\mathbf{\Delta N_t^{in}}$ & $\mathbf{\Delta N_t^{out}}$ \\
\hline
0 & \(N_0\) & \(N_0\) & 0 \\
1 & \(N_0(1 + n_1)\) & \(N_{0}n_1\) & 0 \\
2 & \(N_0(1 + n_1)(1 + n_2)\) & \(N_0(1 + n_1)n_2\) & 0 \\
\vdots & \vdots & \vdots& \vdots \\
\(\mathbf{t}\) &
\(N_0\left( 1 + n_1 \right)\ldots\left( 1 + n_{t} \right)\) \newline
\(=\mathbf{\frac{N_0(1 + n)^t}{ 1 +\frac{N_0}{N}\left\lbrack (1 + n )^t -1\right\rbrack}} \)
\(\mathbf{= :\,N_t^{t < T}}\)
&\(N_0\left( 1 + n_1 \right)\ldots\left( 1 + n_{t - 1} \right)n_t\) \newline
\(=\mathbf{N_{t -1}^{t < T}\,n_t}^{}\)
& \textbf{0} \\
\vdots & \vdots & \vdots& \vdots \\
\(T - 1\) &
\(N_0\left( 1 + n_1 \right)\ldots\left( 1 + n_{T - 1} \right)\) &
\(N_0\left( 1 + n_1 \right)\ldots\left( 1 + n_{T - 2} \right)n_{T - 1}\)
& 0 \\
\hline
\(T\) &
\(N_0\left( 1 + n_1 \right)\ldots\left( 1 + n_T \right) - N_0\) \newline
\(= N_T^{t < T} - N_0\) &
\(N_0\left( 1 + n_1 \right)\ldots\left( 1 + n_{T - 1} \right)n_T\) \newline
\(= N_{T - 1}^{t < T}\ n_T\) & \(\Delta N_0^{in}\) \\
\(T + 1\) &
\(N_0\left( 1 + n_1 \right)\ldots\left( 1 + n_T \right) - N_0\) \newline
\(+ N_0\left( 1 + n_1 \right)\ldots\left( 1 + n_T \right)n_{T + 1}\) \newline
\(- N_0n_1\) \newline
\(= N_0\left( 1 + n_1 \right)\ldots\left( 1 + n_{T + 1} \right)\) \newline
\( - N_0\left( 1 + n_1 \right)\) \newline
\(= N_{T + 1}^{t < T} - N_1^{t < T}\) &
\(N_0\left( 1 + n_1 \right)\ldots\left( 1 + n_T \right)n_{T + 1}\) \newline
\(= N_{T}^{t < T}\ n_{T + 1}\) & \(\mathrm{\Delta}N_{1}^{in}\) \\
\vdots & \vdots & \vdots& \vdots \\
\(\mathbf{t}\) & \(\mathbf{N_t^{t < T}-N_{t - T}^{t < T}}\)
&\(N_0\left( 1 + n_1 \right)\ldots\left( 1 + n_{t - 1} \right)n_t\) \newline
\(\mathbf{=N_{t -1}^{t < T}n_t}\) &\(\mathbf{\Delta N_{t - T}^{in}}\) \\
\hline
\end{tabular}
\normalsize
\caption{Evolution of the number of participants in a system with a quasi-logistically decreasing growth rate $n_t$. The number of members $N_t$ of all investors, the number of new(ly entering) investors $\Delta N_t^{in}$, and the number of exiting investors $\Delta N_t^{out}$ develop iteratively according to equations \eqref{recrelNt} and \eqref{QLNt} to \eqref{QLNtlogfct}. For a derivation of $N_t^{t < T}$, see Note \#\,\ref{fn5}.}
\end{table}

\section{Quasi-logistic capital dynamics\newline in continuous time}\label{AppB}

The discrete-time investment model with quasi-logis\-tic demographic growth as discussed in Sections \ref{Sect2p2} and \ref{Sect3p4p1} (embedding Ponzi scheme-like behaviour) is here transformed into a continuous-time model with an exact solution expressed in terms of the hypergeometric function $_2F_1(\alpha,\beta;\gamma;\delta)$. 

\emph{Demographic evolution}. In discrete time, the number of early participants is calculated by means of the quasi-logistic sequence (starting with $N_0$ and proceeding according to equations \ref{QLNt}-\ref{QLNtlogfct}), and the change of the number of partici\-pants at the beginning of period $t\ge 1$ can be expressed by $\Delta N_t^\text{in}=n_t N_{t-\Delta t} \Delta t$ (equation \ref{QLDNtin2}, with time step $\Delta t=1$), whereas for $t\ge T$ an additional subtractive term $\Delta N_t^\text{out}=\Delta N_{t-T}^\text{in}$ appears.  In the continuous-time limit $\Delta t \rightarrow 0$ (and setting some prefactor $(e^q-1)/(1-N_0/N)\approx q$ for ease of notation), this approximately becomes
\begin{eqnarray}
\dot{N}^\text{in}(t)\,dt &=& qN_0 \,\frac{e^{qt}}{\left( 1+ae^{qt} \right)^2}\,dt \hspace{0.4cm} \forall t  \in \mathbb{R}^{+}\\
\dot{N}^\text{out}(t)\,dt &=& \dot{N}^\text{in}(t-T)\,dt \hspace{1.2cm} t\ge T, 
\end{eqnarray}
with initial number of investors $N_0=N(0)$, initial population growth rate $q=\ln (1+n)$, and a ratio $a\equiv N_0/(N-N_0)$. The number of current investors is $N(t)=N_0+\int_0^t \dot{N}^\text{in}(\tau)\,d\tau$ (for $0\le t < T$) and $N(t)=N_0+\int_0^t \dot{N}^\text{in}(\tau)\,d\tau-\int_T^t \dot{N}^\text{in}(\tau-T)\,d\tau$ (for $t \ge T$), yielding 
\begin{equation}
N(t) =  N_0 \,\begin{cases}
\left( 1+ \frac{1}{a}\left( \frac{1}{1+a} - \frac{1}{1+a e^{qt}} \right) \right)& 0\le t < T\\
\left(      \frac{1}{a}\left( \frac{1}{1+a e^{q(t-T)}} - \frac{1}{1+a e^{qt}} \right) \right) & t \ge T.
\end{cases}
\end{equation}
The function $N(t)$ for $t<T$  is called the \emph{quasi-logistic function}, because in the case $N_0/N \ll 1$ it may be written in the form of the logistic function $N(t)\approx N/\left(1+\frac{N-N_0}{N_0}e^{-qt}\right)$ (see equation \ref{QLNtlogfct}). For $t\ge T$, the graph of the function $N(t)$ changes from sigmoidal to humped. Then the number of investors peaks at time 
\begin{equation}
t_\text{max}=\frac{T}{2}+\frac{1}{q}\ln \left( \frac{N}{N_0} -1\right).
\end{equation}
(This is similar to equation \eqref{QLtNpeak} in the ($\Delta t$$\rightarrow$$0$)-limit.) The inverse function is needed in order to halt the scheme at the time of occurence of a certain number $N^\star$ of remaining investors, say, $N^\star\lesssim 1$; it is found to be 
\begin{equation}
t(N^\star) =\frac{1}{q}\ln \left( \frac{-B+\sqrt{B^2-4AC}}{2A} \right),
\end{equation}
with
\begin{eqnarray}
\nu&=&\frac{N_0}{N},\,\,\,\nu^\star=\frac{N^\star}{N},  \nonumber\\
A&=&\frac{\nu^\star\nu^2}{(1-\nu)^3}e^{-qT}, \,\,\,C=\frac{\nu^\star}{1-\nu}, \nonumber\\
B&=& \frac{\nu}{1-\nu}\left( \left(\frac{\nu^\star}{1-\nu}+1\right) \left( 1+e^{-qT}\right) -2\right)
\end{eqnarray}

\emph{Capital evolution} is governed by the  first-order ordinary differential equation (ODE) 
\begin{equation}
\dot{K}(t)=pK(t)+\dot{C}^\text{in}(t)-\dot{C}^\text{out}(t),
\end{equation}
where a dot represents the first derivative with respect to the time variable $t$, and $p$ is the rate of realized nominal return (formerly denoted by $i$). The cash inflow rate due to deposits is $\dot{C}^\text{in}(t)=I_0 \dot{N}^\text{in}(t)$; the cash outflow rate is $\dot{C}^\text{out}(t)=I_0 \dot{N}^\text{out}(t)+rI_0N(t)$, representing withdrawals and profits, respectively (the latter corresponding to a return on investment with rate $r$). Adopting the integration factor method, the ODE with initial value $K_0=K(0)$ has the following formal solution 
\begin{equation}
K(t) = e^{pt}\left( K_0+\int_0^t \left( \dot{C}^\text{in}(\tau)-\dot{C}^\text{out}(\tau) \right)e^{-pt} d\tau \right).
\end{equation}
Relying on the expressions above and on the indefinite integral
\begin{equation}
\int \frac{e^{ux}}{\left( 1+we^{vx} \right)^n}\,dx = \frac{e^u}{u}
{}_2F_1\left( n, \frac{u}{v};1+\frac{u}{v};-w e^{vx} \right), \label{indefint}
\end{equation}
the explicit solution is found to be
\setlength{\leftmargini}{0cm}
\begin{equation}
\hspace{-1.0cm}
\mathbf{0\le t <T}: \,\,\,K(t) = K_0e^{pt} +D_{<}(t) - P_{<}(t)
\end{equation}
with 
\vspace{-0.2cm}
\begin{eqnarray}
D_{<}(t) &=& \frac{N_0I_0q}{q-p}\left( e^{qt} {}_2F_1\left( 2, 1-\frac{p}{q};2-\frac{p}{q};-a e^{qt} \right) \right.
		\nonumber\\
	   & & \left. -\,{}_2F_1\left( 2, 1-\frac{p}{q};2-\frac{p}{q};-a  \right) \right) \\
P_{<}(t) &=& \frac{N_0I_0 r}{p a} \left( \frac{1+a+a^2}{1+a}  \left( e^{pt} -1 \right) \right. \nonumber \\
	   & & -e^{pt} {}_2F_1\left( 1,-\frac{p}{q};1-\frac{p}{q};-a \right) \nonumber \\
	   & & \left. +\,{}_2F_1\left( 1, -\frac{p}{q};1-\frac{p}{q};-a e^{qt} \right) \right) ,
\end{eqnarray}
representing the aggregated deposits $D_{<}(t)$ (excluding the initial deposits $N_0I_0$ that enter $K_0$) and profits $P_{<}(t)$, respectively.
\begin{equation}
\mathbf{t \ge T}:  \,\,\,K(t) = K_0e^{pt} +D_{>}(t) - P_{>}(t)-W_{>}(t),
\end{equation}
with  
\vspace{-0.2cm}
\begin{eqnarray}
D_{>}(t) &=& D_{<}(t). \\
P_{>}(t) &= &  P_{<}(T) +\frac{N_0I_0 r}{p a} \Bigg[ e^{p(t-T)} \,\times  \nonumber\\
	    & & {}_2F_1\left( 1,-\frac{p}{q};1-\frac{p}{q};-a \right) \nonumber \\
	    & &+\, {}_2F_1\left( 1,-\frac{p}{q};1-\frac{p}{q};-a e^{qt} \right) \nonumber \\
	    & & -\,e^{p(t-T)} {}_2F_1\left( 1,-\frac{p}{q};1-\frac{p}{q};-a e^{qT} \right) \nonumber \\
	    & &-\,  {}_2F_1\left( 1,-\frac{p}{q};1-\frac{p}{q};-a e^{q(t-T)} \right) \Bigg] \nonumber \\
	    & & \\
W_{>}(t) &= & N_0I_0 \Bigg[ 1 +  \frac{q}{q-p} \biggl( e^{q(t-T)} \times \nonumber \\
	    & &   {}_2F_1\left( 2, 1-\frac{p}{q};2-\frac{p}{q};-a e^{q(t-T)} \right)  \nonumber \\
	    & & \left. -\,{}_2F_1\left( 2, 1-\frac{p}{q};2-\frac{p}{q};-a \right) \right)\Bigg].
\end{eqnarray}
These contributing quantities represent the aggregated deposits $D_{>}(t)$, profits $P_{>}(t)$, and withdrawals $W_{>}(t)$, respectively. 

A graphical illustration of the capital evolution of such an investment system in continuous time and under varying conditions is provided in the main body of the paper (Fig.~\ref{Fig6}, dotted lines in the upper three panels); it closely follows the results obtained in discrete time and hence exhibits the characteristic trio of \emph{traffic light scenarios}.

\newpage 
\section{Table of Contents}\label{AppC}
\renewcommand{\contentsname}{}
\tableofcontents

\end{document}